\documentclass[pdflatex,sn-mathphys-num]{sn-jnl}

\usepackage{graphicx}%
\usepackage{multirow}%
\usepackage{amsmath,amssymb,amsfonts}%
\usepackage{amsthm}%
\usepackage{mathrsfs}%
\usepackage[title]{appendix}%
\usepackage{xcolor}%
\usepackage{textcomp}%
\usepackage{manyfoot}%
\usepackage{booktabs}%
\usepackage{algorithm}%
\usepackage{algorithmicx}%
\usepackage{algpseudocode}%
\usepackage{listings}%
\usepackage{subfigure}
\usepackage{longtable} 
\usepackage{pdflscape}
\usepackage{rotating}
\usepackage{todonotes}

\begin{document}

\title[Heterogeneous Min-Max TSP]{Heterogeneous Min-Max Multi-Vehicle Multi-Depot Traveling Salesman Problem: Heuristics and Computational Results}


\author*[1]{\fnm{Deepak Prakash} \sur{Kumar}}\email{deepakprakash1997@gmail.com}

\author[2]{\fnm{Sivakumar} \sur{Rathinam}}\email{srathinam@tamu.edu}

\author[3]{\fnm{Swaroop} \sur{Darbha}}\email{dswaroop@tamu.edu}
\author[4]{\fnm{Trevor}\sur{Bihl}}\email{trevor.bihl.2@afrl.af.mil}

\affil*[1, 2, 3]{Department of Mechanical Engineering, Texas A\&M University, 3123 TAMU, College Station, Texas 77843, USA}
\affil[4]{Air Force Research Laboratory, Sensing Management Branch, Dayton, OH 45433, USA}


\abstract{In this paper, a heuristic for a heterogeneous min-max multi-vehicle multi-depot traveling salesman problem is proposed, wherein heterogeneous vehicles start from given depot locations and need to cover a given set of targets. In the considered problem, vehicles can be structurally heterogeneous due to different vehicle speeds and/or functionally heterogeneous due to different vehicle-target assignments originating from different sensing capabilities of vehicles. The proposed heuristic for the considered problem has three stages: an initialization stage to generate an initial feasible solution, a local search stage to improve the incumbent solution by searching through different neighborhoods, and a perturbation/shaking stage, wherein the incumbent solution is perturbed to break from a local minimum. In this study, three types of neighborhood searches are employed. Furthermore, two different methods for constructing the initial feasible solution are considered, and multiple variations in the neighborhoods considered are explored in this study. The considered variations and construction methods are evaluated on a total of 128 instances generated with varying vehicle-to-target ratios, distribution for generating the targets, and vehicle-target assignment and are benchmarked against the best-known heuristic for this problem. Two heuristics were finally proposed based on the importance provided to objective value or computation time through extensive computational studies.}


\maketitle

\section{Introduction} \label{sec1}

The Traveling Salesman Problem (TSP) \cite{TSP_book_applegate} and its variants are some of the popular vehicle routing problems addressed due to their practical applications. In the traveling salesman problem, a salesman/vehicle needs to cover a given set of targets in a graph such that each target is visited exactly once, the salesman returns to the starting location at the end of the tour, and the tour cost is minimal. Various variants of this problem exist in the literature \cite{TSP_variants}. One of the variants of interest is the multi-vehicle Traveling Salesman Problem (TSP), wherein multiple vehicles are deployed from a given set of depots (that are not necessarily distinct). The multi-vehicle TSP variant has substantial civilian applications, such as medicine and package delivery, and military applications, such as surveillance. In practice, when a set of vehicles are considered to be utilized for a particular application, they do not always necessarily need to be of the same make or design. In this regard, it is necessary to consider routing a heterogeneous set of vehicles to accomplish the underlying mission. Typically, two types of heterogeneity are considered in vehicles: functional and structural heterogeneity. Functional heterogeneity in vehicles arises due to a different set of sensors that are mounted on vehicles, which in turn leads to a particular set of targets in the graph to be necessarily visited by a particular vehicle. On the other hand, structural heterogeneity in vehicles could arise from the vehicles traveling at different speeds or having different kinematic constraints, such as the turning radius in the case of a Dubins vehicle \cite{Dubins}.

For the multi-vehicle TSP, two variants that differ in the objective function are typically considered: a min-sum variant and a min-max variant. In the min-sum variant, tours must be constructed for the vehicles such that the total cost of tours for the vehicles is minimized. The min-sum variant has been studied extensively, and efficient algorithms for the same exist. For example, the optimal solution for instances of the min-sum variant can be obtained by solving a formulated mixed-integer linear program (MILP) in \cite{algorithms_heterogeneous_multi_depot_multi_vehicle}. The authors in \cite{algorithms_heterogeneous_multi_depot_multi_vehicle} showed that instances with $100$ targets and $5$ vehicles can be solved to optimality using branch and cut in $300$ seconds on average, due to the strong lower bounds that can be obtained from the linear program relaxation. Alternatively, the problem can be solved to near optimality using the state-of-the-art heuristic for a single-vehicle TSP, known as LKH \cite{LKH_TSP_solver}, by performing a graph transformation, known as the Noon-Bean transformation \cite{noon_bean}. Such a transformation has been demonstrated for the min-sum problem for heterogeneous vehicles in \cite{todays_TSP}, wherein the min-sum multi-vehicle problem gets transformed to an equivalent asymmetric TSP, which can be solved using LKH. In \cite{todays_TSP}, the authors showed that instances with about $40$ targets and $10$ heterogeneous vehicles modeled as Dubins vehicles due to turning constraints could be solved within $50$ seconds to within about $10\%$ of the optimal solution.

On the other hand, the min-max variant is a more difficult variant to tackle due to poor lower bounds obtained from linear program (LP)-based relaxations of the corresponding MILPs. Typically, this variant is studied under two variations: a single-depot variant and a multi-depot variant. In the single-depot variant, all vehicles must start and end at a given depot, and many studies have considered this variant of the problem \cite{survey_mtsp, Memetic_search_min_max_TSP}. For the single-depot variant, approximation algorithms exist in the literature, which typically originates from the algorithm given in \cite{Frederickson} for the problem with min-max problem with homogeneous vehicles. The authors proposed a $1 + F - \frac{1}{k}$-factor algorithm, where $k$ denotes the number of vehicles and $F$ denotes the approximation ratio for the problem with a single vehicle. Using the well-known Christofides algorithm, which yields a $\frac{3}{2}$-approximation algorithm for the single vehicle TSP \cite{Christofides}, a $\frac{5}{2} - \frac{1}{k}$ approximation ratio algorithm is obtained for the min-max problem with homogeneous vehicles starting at a single depot. Studies such as \cite{Battistini} and \cite{Prasad} utilize this approximation-ratio algorithm and extend it for heterogeneous variants of the problem, wherein the vehicles are heterogeneous due to different task costs in the former, whereas the vehicles are functionally heterogeneous in the latter.

In contrast, the multi-depot variant of the min-max problem is typically tackled using heuristics due to the existence of poor lower bounds. Feasible solutions for the min-max multi-vehicle multi-depot TSP with homogeneous vehicles were generated using four different heuristics in \cite{min_max_vrp_LB_based_load_balancing}. The authors observed that a heuristic based on an LP-based load balancing method and a heuristic based on region partitioning performed the best for $13$ instances. In \cite{ant_colony_min_max}, the authors proposed an ant-colony-based method, benchmarked it against the LP-based heuristic from \cite{min_max_vrp_LB_based_load_balancing}, and observed that the heuristic performed better than the LP-based method for a few instances wherein the targets were generated from a uniform distribution. In \cite{MD_algorithm}, the authors proposed a heuristic, termed as the ``MD" algorithm, which is a three-staged algorithm. In the first stage, the authors utilize the LP-based load balancing method from \cite{min_max_vrp_LB_based_load_balancing} to obtain a quick initial feasible allocation of targets to vehicles and utilize LKH to obtain the tours for each vehicle. A local search is then performed on a neighborhood to improve the solution based on a search in a neighborhood obtained by switching a target from the vehicle with the highest tour cost to another vehicle.
A perturbation step is then used to break from a local minimum. The proposed algorithm was benchmarked against three other algorithms, including a Variable Neighborhood Search algorithm \cite{Mladenovic1997}, which is a traditional heuristic used for producing high-quality feasible solutions. A memetic algorithm was recently proposed in \cite{Memetic_search_min_max_TSP}, in which multiple initial solutions are generated, and neighborhood searches are performed to improve these initial solutions. The authors were able to produce the best-known results for a majority of instances generated in \cite{MD_algorithm}.

Unlike the min-max multi-vehicle multi-depot TSP with homogeneous vehicles, fewer studies have explored the heterogeneous variant of the problem. In \cite{3_approx_algo_2_vehicle_TSP}, the authors propose a $\frac{3k}{2}$-approximation ratio algorithm, wherein the bound for the feasible solution obtained increases with the number of vehicles $k$. The algorithm is currently the best-known approximation ratio algorithm for such a problem. A heuristic based on a primal-dual algorithm was proposed in \cite{heuristic_struct_het_distict_dep} for vehicles with different speeds. However, the authors assume that the vehicles have distinct speeds and start at distinct depots. Recently, a heuristic with three stages was proposed in \cite{AIAA_modified_MD}, wherein the MD algorithm proposed for the homogeneous problem in \cite{MD_algorithm} was generalized for the heterogeneous problem. In this study, the authors considered vehicles to be functionally heterogeneous due to different vehicle-target assignments and structurally heterogeneous due to different (and not necessarily distinct) vehicle speeds. The generalized heuristic was evaluated on instances that were randomly generated for three vehicles and thirty targets. However, in the local search stage of the heuristic, only one neighborhood search was employed to improve the incumbent solution. Furthermore, in the heuristic, LKH \cite{LKH_TSP_solver} was used to optimize the tours of the vehicles for every explored solution in the neighborhood. However, such an implementation is computationally expensive (as will be observed from the results presented in Section~\ref{sect: results}), and no variations in the construction heuristic were explored. To this end, the main contributions of this paper are as follows:
\begin{enumerate}
    \item A three-staged heuristic inspired by the MD algorithm \cite{MD_algorithm} is proposed, wherein three diverse neighborhoods are utilized in the local search stage (the second stage) to obtain high-quality feasible solutions for the heterogeneous min-max multi-vehicle multi-depot TSP.
    \item A set of $128$ instances are generated with varying vehicle-target ratios, vehicle-target allocations, and distributions from which the targets are generated, with instances considering as large as $500$ targets and $20$ vehicles.
    \item Two types of construction methods for the initial solution, and $18$ variations in the first two neighborhoods and $16$ variations in the third neighborhood are considered and are studied through comprehensive numerical simulations. Furthermore, the results are benchmarked against the heuristic in \cite{AIAA_modified_MD}, which is the current best-known heuristic for the studied problem.
    \item Two heuristics are proposed based on computational studies performed on the developed instances, one of which yields solutions with a better objective value, whereas the other heuristic is computationally efficient. In particular, the heuristic selected for the best objective value produces a better solution or equivalent solution compared to the benchmarking heuristic on $71$ and $41$ instances, respectively, out of $128$ instances. On the other hand, the heuristic selected based on the computation time produces a feasible solution within $10$ minutes on $127$ out of the $128$ instances, and produced a solution with a better or equal objective value on $64$ and $40$ instances, respectively, out of $128$ instances compared to the benchmarking heuristic.
\end{enumerate}

\section{Problem Description and Notation}

Consider a set of targets $T$ that need to be covered by $k$ vehicles, which start at depots $d_1, d_2, \cdots, d_k$, with speeds $v_1, v_2, \cdots, v_k$. It should be noted that the depots and the vehicle speeds are not necessarily distinct. Further, let set $R_i$ denote the set of targets in $T$ that need to be necessarily covered by vehicle $i$ for $i \in \{1, 2, \cdots, k \}$. The sets $R_1, \cdots, R_k$ are considered to be mutually disjoint, i.e., $R_i \cap R_j = \emptyset$ for all $i, j = 1, 2, \cdots, k,$ and $i \neq j$. In this study, it is desired to construct tours for the $k$ vehicles such that
\begin{itemize}
    \item The $i\textsuperscript{th}$ vehicle's tour starts and ends at depot $d_i$ for all $i = 1, 2, \cdots, k,$
    \item The targets in set $R_i$ are necessarily covered by the $i\textsuperscript{th}$ vehicle for $i = 1, 2, \cdots, k$,
    \item Each target in the set $T \setminus \{R_1 \cup R_2 \cup \cdots \cup R_k \}$ are covered by exactly one vehicle, and
    \item The maximum tour time of the $k$ vehicles is minimized.
\end{itemize}
Let $V_i = T \cup \{d_i\}$ denote the set of vertices that can be traversed by the $i\textsuperscript{th}$ vehicle and $E_i$ denote the set of edges connecting vertices in $V_i$. In this study, it is assumed that the cost of edge $(m, n) \in E_i$ equals the time taken to travel, which is the Euclidean distance between $m$ and $n$ divided by the vehicle speed $v_i$, for all $i \in \{1, 2, \cdots, k\}$. Further, the graph $G_i$ associated with the $i\textsuperscript{th}$ vehicle is assumed to be undirected and complete. The MILP formulation for the considered problem is provided in \cite{3_approx_algo_2_vehicle_TSP}.

\textbf{Remark:} Though the cost of an edge is considered to be the time taken by the vehicle to traverse it, the proposed heuristic in this paper and the MILP formulation in \cite{3_approx_algo_2_vehicle_TSP} are applicable for graphs that are symmetric and whose edge costs satisfy the triangle inequality.

\section{Heuristic}

\begin{figure}[htb!]
    \centering
    \includegraphics[width = 0.6\linewidth]{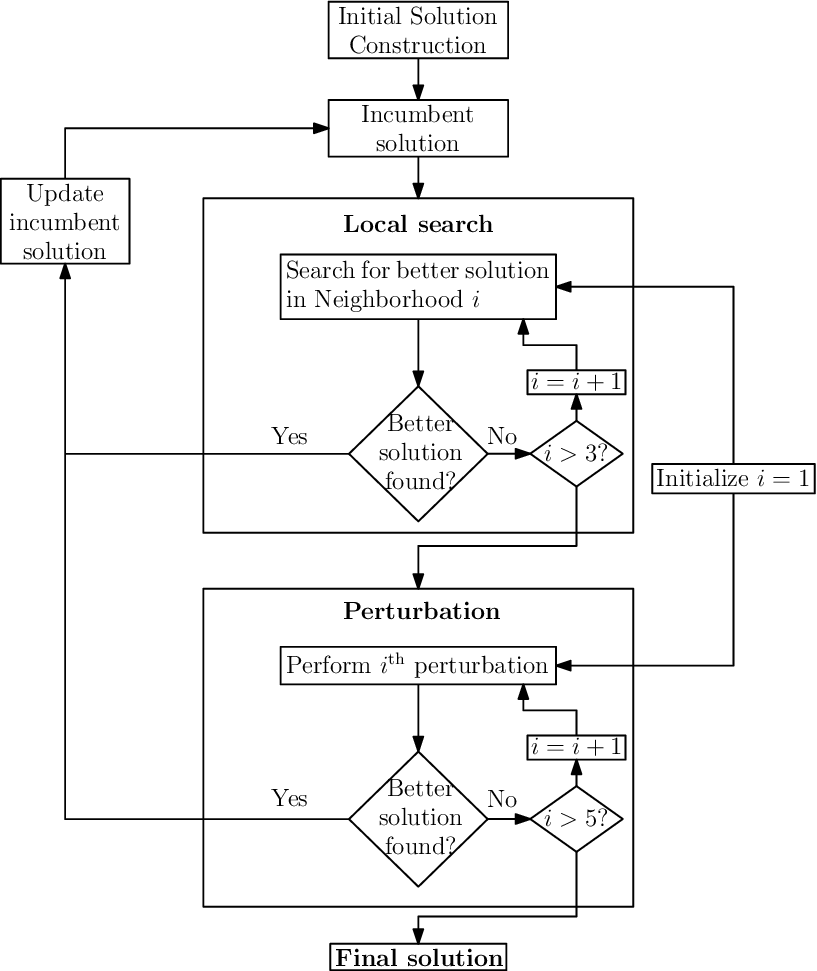}
    \caption{Overview of stages of the proposed heuristic}
    \label{fig: overview_heuristic}
\end{figure}

The proposed heuristic in this study is inspired by the MD algorithm \cite{MD_algorithm} and the modified MD algorithm for the considered problem \cite{AIAA_modified_MD}, and consists of three stages:
\begin{itemize}
    \item Initialization: An initial feasible solution is generated.
    \item Local search: The incumbent solution is improved using neighborhoods.
    \item Perturbation: The incumbent solution is perturbed at random to break from a local minimum.
\end{itemize}
A depiction of these three stages to improve the solution is shown in Fig.~\ref{fig: overview_heuristic}. The description of these three stages follows.

\subsection{Initialization} \label{subsect: initialisation_step}

In the initialization stage, an initial feasible solution for the considered problem is obtained. In this paper, two approaches are considered for constructing an initial feasible solution. In the first approach, a linear program-based load balancing is used, which was proposed in \cite{AIAA_modified_MD}, and can be solved with solvers such as Gurobi \cite{gurobi}. In this approach, it is desired to allocate a desired number of targets based on the vehicle speed so that the cost of assigning targets to the vehicles is minimized. While such an algorithm was observed to perform well in \cite{AIAA_modified_MD} for instances without vehicle-target assignment, the initial solution obtained for instances with vehicle-target assignment had a larger deviation with respect to the optimum for small instances. In this regard, an additional construction heuristic is considered in this study, which is inspired by heuristics studied in \cite{dubins_tsp_heuristics} for a Traveling Salesman Problem for a Dubins vehicle \cite{Dubins}\footnotemark. The construction heuristic proposed is based on recursive insertion, wherein the vehicle tours are constructed step by step by inserting one target at a time.

\footnotetext{A Dubins vehicle is a vehicle that moves forward at a unit speed and has a bound on the rate of change of the heading angle.}

\begin{figure}[htb!]
    \centering
    \subfigure[Initial tours generated using pre-assigned targets]{\includegraphics[width = 0.48\textwidth]{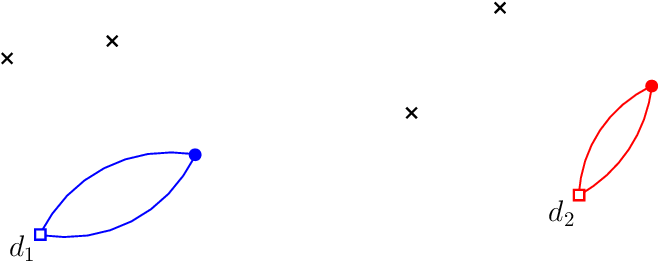}\label{subfig: construction_before}} \hfill
    \subfigure[Tours of vehicles after first step of recursive insertion]{\includegraphics[width = 0.48\textwidth]{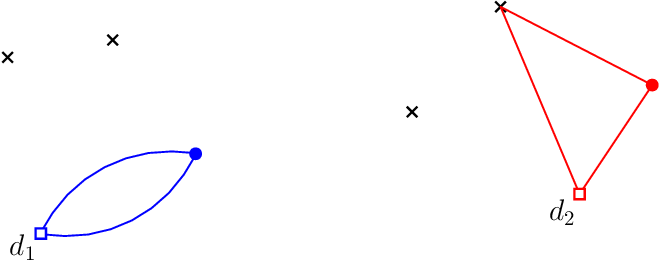}\label{subfig: construction_after}}
    \caption{Depiction of recursive insertion method for constructing an initial feasible solution}
    \label{fig: recursive_insertion_initial_solution}
\end{figure}

In the first step of the recursive insertion method, the targets that are pre-assigned to each vehicle, i.e., in $R_i$ for $i = 1, 2, \cdots, k,$ are used to construct the initial tour for each vehicle using LKH, as depicted in Fig.~\ref{subfig: construction_before}. Following this initial construction of tours, the following sequence of steps is repeated till there are no more targets in $T \setminus \{R_1 \bigcup R_2 \bigcup \cdots \bigcup R_k \}$ to be assigned:
\begin{enumerate}
    \item The vehicle with the least tour cost is selected.
    \item The cost of insertion is computed for each target in $T \setminus \{R_1 \bigcup R_2 \bigcup \cdots \bigcup R_k \}$ that has not been assigned to a vehicle, whose definition follows. Suppose target $t$ is attempted to be inserted in the $j\textsuperscript{th}$ vehicle, whose tour is initially given by $(d_j, u_1^j, u_2^j, \cdots, u_l^j, u_{l + 1}^j, \cdots, u_{n_j - 1}^j, d_j).$ Let target $t$ be attempted to be inserted between $u_l^j$ and $u_{l + 1}^j$. The insertion cost is defined, similar to the definition in \cite{AIAA_modified_MD}, as
    \begin{align} \label{eq: insertion_cost_single_vehicle_definition}
        \text{insertion cost}_{t, j, (u_l^j, u_{l + 1}^j)} = \frac{dist(u_l^j, t)}{v_j} + \frac{dist(t, u_{l + 1}^j)}{v_j} - \frac{dist(u_l^j, u_{l + 1}^j)}{v_j},
    \end{align}
    where $dist$ denotes the Euclidean distance between two considered targets. Since target $t$ can be inserted between any pair of vertices in vehicle $j$'s tour, the insertion cost for target $t$ is defined to be the minimum among all such costs.

    Among all targets that have not been assigned in $T \setminus \{R_1 \bigcup R_2 \bigcup \cdots \bigcup R_k \}$, the target with the least insertion cost is selected and is inserted in the tour of the vehicle with the least tour cost. The location of insertion in the chosen vehicle's tour corresponds to the location that yields the least insertion cost using Eq.~\eqref{eq: insertion_cost_single_vehicle_definition}. A depiction of the modified tours of the vehicles after this step is shown in Fig.~\ref{subfig: construction_after}.
    \item The considered target for insertion is then removed from the common set of targets $T \setminus \{R_1 \bigcup R_2 \bigcup \cdots \bigcup R_k \}$ since it has been assigned.
\end{enumerate}
Finally, once all targets in $T \setminus \{R_1 \bigcup R_2 \bigcup \cdots R_k \}$ have been assigned to a vehicle, the LKH heuristic is used to optimize all vehicle tours. In this regard, an initial feasible solution is obtained for the considered problem.

\subsection{Local search}

A local search is the second stage of the heuristic and is performed to improve the incumbent solution. In particular, in the local search, the tour of the vehicle with the highest tour time, termed the maximal vehicle, is attempted to be improved. For this purpose, three neighborhoods are used, which are
\begin{itemize}
    \item Neighborhood 1: Target switch. In this neighborhood, a target is attempted to be removed from the maximal vehicle and inserted into another vehicle.
    \item Neighborhood 2: Target swap. In this neighborhood, a target from the maximal vehicle is attempted to be swapped with a target from another vehicle.
    \item Neighborhood 3: Multi-target swap. In this neighborhood, a set of consecutive targets from the maximal vehicle is attempted to be swapped with a set of consecutive targets from another vehicle.
\end{itemize}
The description of these neighborhoods and the variations considered for these neighborhoods follows.

\subsubsection{Neighborhood 1: target switch}

\begin{figure}[htb!]
    \centering
    \subfigure[Vehicle tours before removing target $t$]{\includegraphics[width = 0.32\textwidth]{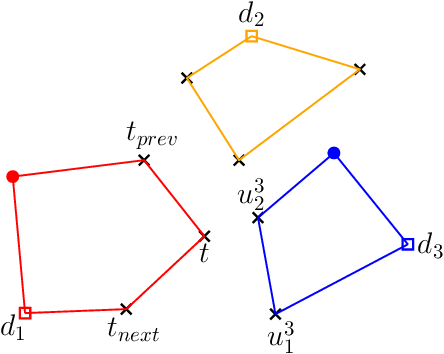}} \hfill
    \subfigure[Maximal vehicle's feasible tour after removing target $t$]{\includegraphics[width = 0.32\textwidth]{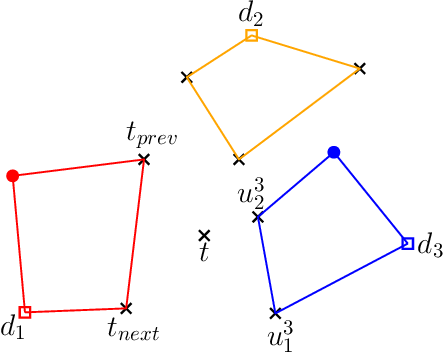}} \hfill
    \subfigure[Vehicle tours after inserting target $t$ in ``blue" vehicle between $p$ and $p'$]{\includegraphics[width = 0.32\textwidth]{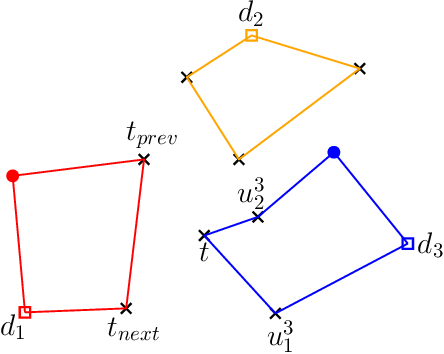}}
    \caption{Illustration of neighborhood 1. Here, the ``red" vehicle is the maximal vehicle, targets colored in black are common targets, and shaded discs represent targets with vehicle-target assignments.}
    \label{fig: neighborhood_1}
\end{figure}

In the target switch neighborhood, a target is attempted to be removed from the maximal vehicle (denoted with index $i$) and inserted into another vehicle, as shown in Fig.~\ref{fig: neighborhood_1}. To this end, it is desired to first order the targets in the maximal vehicle to speed up the search process. Since it is desired to reduce the tour time for the maximal vehicle, a ``savings" metric is introduced, similar to \cite{AIAA_modified_MD}. This metric estimates the change in the tour time with the removal of a particular target from the considered vehicle. The savings associated with removing target $t$ from vehicle $i$, the maximal vehicle, is defined as
\begin{align} \label{eq: definition_savings_metric}
    \text{savings}_{t, i} = \frac{dist(t_{prev}, t)}{v_i} + \frac{dist(t, t_{next})}{v_i} - \frac{dist(t_{prev}, t_{next})}{v_i},
\end{align}
where $t_{prev}$ and $t_{next}$ are vertices covered before and after target $t$, respectively, by vehicle $i$, and $v_i$ is the speed of the $i\textsuperscript{th}$ vehicle.
In Fig.~\ref{fig: neighborhood_1}, $\text{savings}_{t, i}$ denotes the time saved by traveling from $t_{prev}$ to $t_{next}$ directly instead of going through target $t$. All targets in the maximal vehicle, except for those targets in $R_i$, are sorted in the decreasing order of the savings metric since it is desired to remove ``costly" targets for the maximal vehicle. 

Having sorted the targets in the maximal vehicle, the target with the highest savings, denoted by $t$, is first attempted to be removed. The removed target must be covered by one of the other vehicles since all targets in the graph must necessarily be covered by a vehicle. In this regard, a decision needs to be made on the priority of considering vehicles to visit target $t$. Since the best method of picking a vehicle is not immediately apparent, three choices/metrics have been considered for sorting the vehicles (except for the maximal vehicle) to cover target $t$:
\begin{enumerate}
    \item Least actual tour cost: In this metric, vehicles are sorted based on increasing tour time.
    \item Least insertion cost: In this metric, vehicles are sorted based on increasing insertion cost. It should be noted that a similar metric has been used in \cite{MD_algorithm} for the homogeneous problem and in \cite{AIAA_modified_MD} for the heterogeneous problem. 
    \item Least estimated tour cost: In this metric, vehicles are sorted based on increasing estimated tour cost (time) after inserting target $t$.
\end{enumerate}
While sorting the vehicles based on the first metric, i.e., increasing tour time, is immediate, since the tour costs of all vehicles in the current solution are available, the insertion cost, defined in Section~\ref{subsect: initialisation_step}, is utilized for sorting the vehicles (except for the maximal vehicle) through the other two metrics. It should be noted that using the insertion cost computed for inserting target $t$ into vehicle $j$, vehicle $j$'s estimated tour time after insertion can be computed to be
\begin{align} \label{eq: estimated_tour_cost_insertion_N1}
    \text{estimated tour cost}_{j, t} &= \text{previous tour cost}_j + \text{insertion cost}_{t, j}.
\end{align}
In this regard, all vehicles, except the maximal vehicle, can be sorted using the third metric choice as well.

Having sorted the vehicles based on the chosen metric corresponding to the selected target $t$ from the maximal vehicle, the steps considered to determine if a better solution is obtained or not are as follows:
\begin{enumerate}
    \item The first vehicle in the sorted list is initially considered for insertion. Suppose vehicle $j$ is considered for insertion.
    \item The estimated tour cost for the maximal vehicle ($i$) is computed using the savings associated with target $t$ as
    \begin{align}
        \text{estimated tour cost}_{i, t} &= \text{previous tour cost}_i - \text{savings}_{t, i}.
    \end{align}
    Furthermore, the estimated tour cost for vehicle $j$ is obtained using Eq.~\eqref{eq: estimated_tour_cost_insertion_N1}.
    \item Using the new estimated tour costs for vehicle $i$ and $j$, the maximum estimated tour time associated with the new solution is compared with the maximum tour time in the incumbent solution. If the obtained solution corresponds to a lower maximum tour time, then the tours of vehicles $i$ and $j$ are optimized using the LKH heuristic \cite{LKH_TSP_solver}. This is because if a better solution is obtained based on the proxy costs, which are the savings and insertion costs, then a better solution will definitely be obtained using LKH since it yields a near-optimal tour for the two vehicles. If a better solution is obtained from this neighborhood, the incumbent solution is updated, and the search in this neighborhood is stopped.
    
    On the other hand, if a better solution is not obtained using the estimated tour costs, the next vehicle in the sorted list of vehicles based on the chosen metric is considered for inserting target $t$, and the same steps are repeated.
\end{enumerate}

\textbf{Remark:} To speed up the search in this neighborhood, a parameter $n$ is considered, which considers only the top $n$ vehicles in the sorted list for insertion.

\begin{figure}[htb!]
    \centering
    \includegraphics[width = 0.55\linewidth]{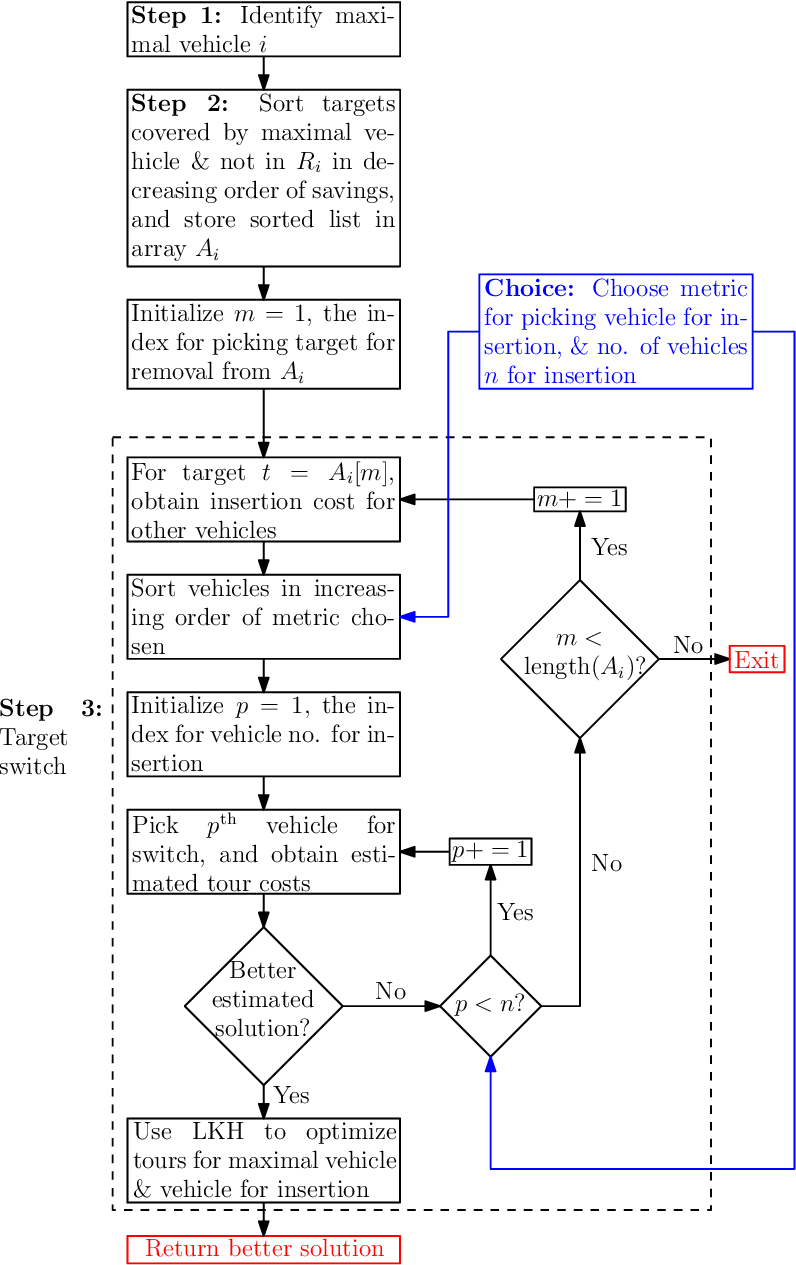}
    \caption{Overview of target switch neighborhood with variations}
    \label{fig: switch_neighborhood}
\end{figure}

If a better solution is not obtained for inserting target $t$ into another vehicle, then the next target in the sorted list of targets in the maximal vehicle ($i$) based on the savings metric is considered. The same set of steps described above are then performed to insert the considered target into another vehicle. The search in this neighborhood is continued till a better solution based on the estimated tour costs is obtained or all targets in the maximal vehicle have been considered. An overview of the implementation of the described neighborhood is shown in Fig.~\ref{fig: switch_neighborhood}.


\textbf{Remark:} It should be noted that while this switch neighborhood is similar to the local search neighborhood in \cite{AIAA_modified_MD}, there are three major differences:
\begin{enumerate}
    \item Two variations in the metric choice for choosing the vehicle for insertion are considered in addition to the least insertion cost metric used in \cite{AIAA_modified_MD}, since the best metric choice for insertion is not immediately apparent.
    \item The number of vehicles for insertion can be varied, whereas in \cite{AIAA_modified_MD}, only one vehicle was considered for insertion. It will later be shown through numerical simulations in Section~\ref{sect: results} that the number of vehicles chosen impacts the improvement in the objective value.
    \item The estimated tour costs are used to check if a better solution is obtained, as opposed to recomputing the vehicle tours every time using LKH, as done in \cite{AIAA_modified_MD}. This method provides significant computational benefits, which will be revisited based on computational results later in this paper.
\end{enumerate}

\subsubsection{Neighborhood 2: target swap}

\begin{figure}[htb!]
    \centering
    \subfigure[Vehicle tours before removing target $t$]{\includegraphics[width = 0.32\textwidth]{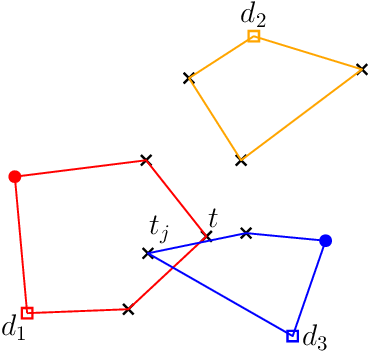}} \hfill
    \subfigure[Maximal vehicle's feasible tour after inserting target $t$ in ``blue" vehicle \label{subfig: vehicle_tours_nbhd_2}]{\includegraphics[width = 0.32\textwidth]{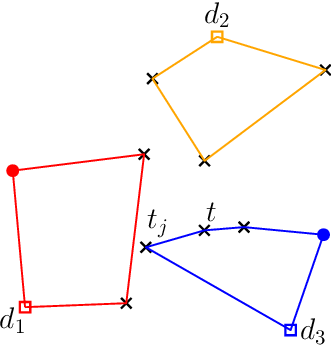}} \hfill
    \subfigure[Vehicle tours after inserting target $t_j$ from ``blue" vehicle into maximal vehicle]{\includegraphics[width = 0.32\textwidth]{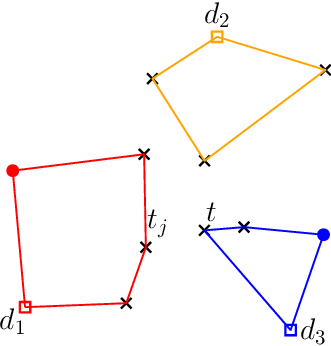}}
    \caption{Illustration of neighborhood 2}
    \label{fig: neighborhood_2}
\end{figure}

In the target swap neighborhood, a target from the maximal vehicle $i$ is attempted to be swapped with a target from another vehicle. The idea behind this neighborhood is to replace a costly target for vehicle $i$ with a cheaper target such that the maximum tour time is reduced, as depicted in Fig.~\ref{fig: neighborhood_2}. To this end, similar to the previous neighborhood, the targets covered by vehicle $i$, except for targets in $R_i$, are first sorted in the decreasing order of the savings metric given in Eq.~\eqref{eq: definition_savings_metric}. First, the target with the highest savings, denoted with $t$, with savings denoted by savings$_{t, i}$ is picked from the maximal vehicle and is attempted to be swapped with another vehicle. Similar to the previous neighborhood, the vehicles to be considered for the target swap, except for vehicle $i,$ are sorted based on one of three metric choices:
\begin{enumerate}
    \item Increasing tour time, which corresponds to the least actual tour cost metric.
    \item Increasing insertion cost for inserting target $t$, which corresponds to the least insertion cost metric, or
    \item Increasing estimated tour time after inserting target $t$, which corresponds to the least estimated tour cost metric.
\end{enumerate}
Similar to the previous neighborhood, the insertion cost for all vehicles defined in Section~\ref{subsect: initialisation_step} is utilized to sort the vehicles based on the second and third metric choices.

\textbf{Remark:} A minor difference in the vehicle choice in this neighborhood compared to the target switch neighborhood is that a vehicle $j$ can be picked for target swap only if it covers at least one target in the current solution that can be covered by another vehicle, i.e., at least one target in $T\setminus R_j$ is covered by vehicle $j$.

\textbf{Remark:} Similar to the previous neighborhood, a parameter $n$ is considered for the neighborhood to speed up the search. The parameter $n$ limits the number of vehicles to be considered for swapping with a selected target from the maximal vehicle.

\begin{figure}[htb!]
    \centering
    \includegraphics[width = 0.6\linewidth]{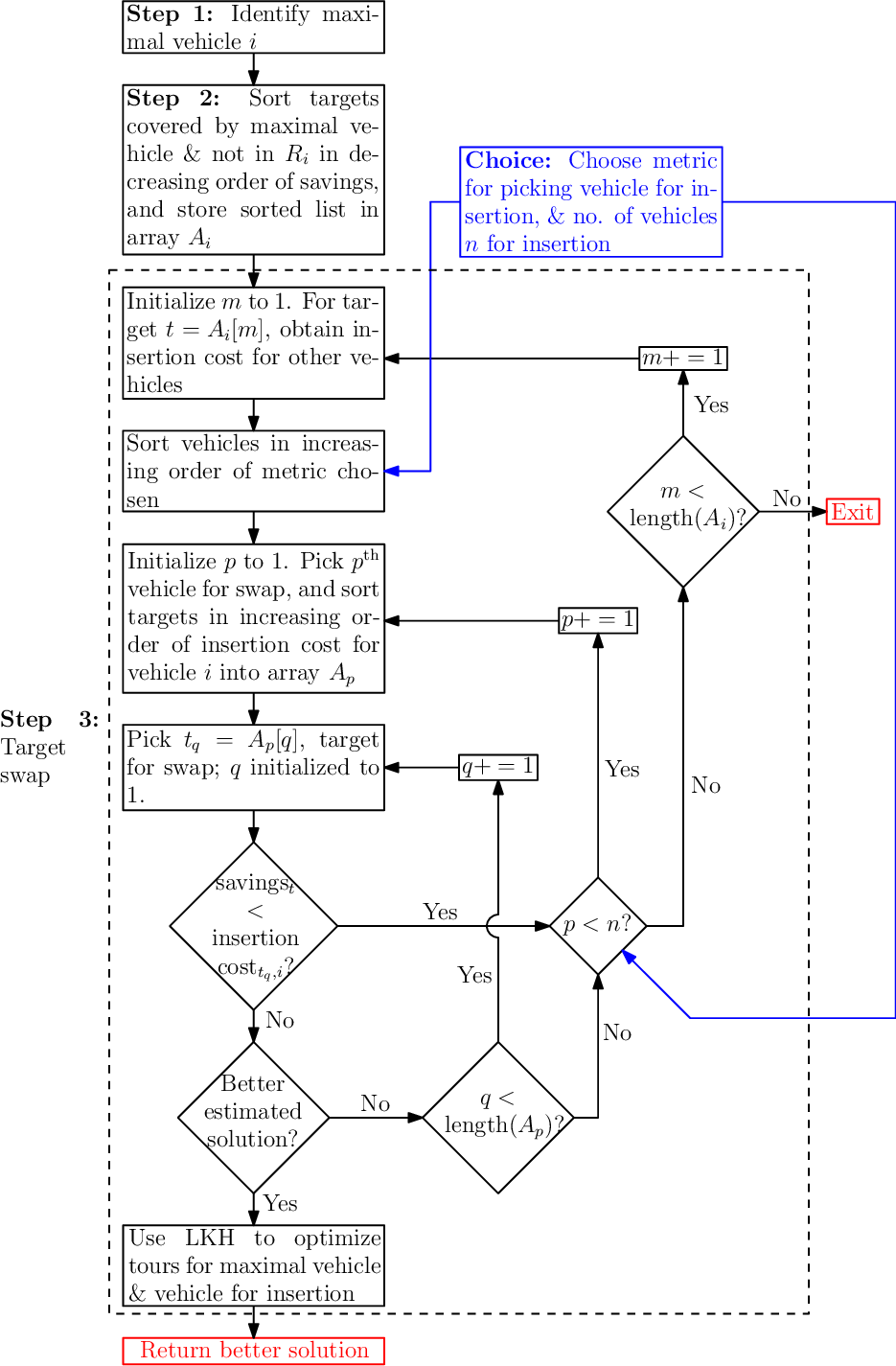}
    \caption{Overview of swap neighborhood with variation}
    \label{fig: swap_neighborhood}
\end{figure}

Having sorted the vehicles based on the chosen metric for attempting to swap with target $t$ from the maximal vehicle, the following steps are then performed to determine if a better solution can be obtained:
\begin{enumerate}
    \item The first vehicle in the sorted list is initially considered for the target swap. Let vehicle $j$ denote the considered vehicle.
    \item Target $t$ is then removed from vehicle $i$'s tour and inserted in the location that yields the least insertion cost, denoted by ``insertion cost$_{t, j}$", based on Eq.~\eqref{eq: insertion_cost_single_vehicle_definition}, as shown in Fig.~\ref{subfig: vehicle_tours_nbhd_2}.
    \item It is now desired to remove a target from vehicle $j$'s tour and insert it in vehicle $i$'s tour. Since a ``cheap" target for vehicle $i$ is desired to be inserted, all targets in vehicle $j$, except for $t$ and targets in $R_j$, are sorted in the increasing order of insertion cost for vehicle $i$. It should be noted that this insertion cost is based on the feasible tour constructed for vehicle $i$ after removing target $t$, as shown in Fig.~\ref{subfig: vehicle_tours_nbhd_2}. Suppose $t_j$ is the target corresponding to the least insertion cost, denoted by ``insertion cost$_{t_j, i}$". The estimated tour costs of the maximal vehicle and vehicle $j$ corresponding to swapping target $t$ with the considered target $t_j$ can be computed to be
    \begin{align*}
        \text{estimated tour cost}_i &= \text{previous tour cost}_i - \text{savings}_{t, i} + \text{insertion cost}_{t_j, i}, \\
        \text{estimated tour cost}_j &= \text{previous tour cost}_j + \text{insertion cost}_{t, j} - \text{savings}_{t_j, j},
    \end{align*}
    where savings$_{t_j, j}$ denotes the savings obtained associated with removing target $t_j$ from vehicle $j$, whose value can be computed using Eq.~\eqref{eq: definition_savings_metric}.
    The following steps are then checked to determine if a better solution is obtained, another target needs to be considered for a swap, or another vehicle must be explored:
    \begin{enumerate}
        \item 
        If $\text{savings}_{t, i}$ is less than $\text{insertion cost}_{t_j, i},$ then the considered target swap will not yield a better solution based on the proxy costs since the estimated tour cost for vehicle $i$ will be higher than the incumbent solution's cost. In this case, other targets in vehicle $j$ are not explored since they have been sorted in the increasing order of the insertion cost for vehicle $i$. Hence, the next vehicle is considered for the swap with target $t$ if not more than $n$ vehicles have been explored for the swap corresponding to target $t$, and Step~1 is restarted with the newly considered vehicle. However, if more than $n$ vehicles have been explored, then the next target in the sorted list based on savings for vehicle $i$ is considered, and all vehicles are then sorted based on the chosen metric. Step~1 is then restarted for the new target.
        \item If $\text{savings}_{t, i}$ is greater than or equal to $\text{insertion cost}_{t_j, i},$ then the tour cost of the maximal vehicle has reduced or at least remains the same based on the proxy costs. Now, the objective value for the considered solution is compared with the incumbent solution's cost. If the objective value, i.e., maximum tour time, is less than the incumbent cost, then LKH is used to optimize the tours of vehicle $i$ and $j$, and the incumbent solution is updated. Neighborhood $1$ is then explored to improve the solution once again (refer to Fig.~\ref{fig: overview_heuristic}).
    
        If a better solution is not obtained based on the proxy costs, then the next target in vehicle $j$'s tour based on the insertion cost is considered for swapping and Step~3 is repeated. If no more targets are available to be explored in vehicle $j$'s tour, then the next vehicle is considered for the swap if no more than $n$ vehicles have been explored for the swap corresponding to target $t$. Using the newly considered vehicle, Step~1 is restarted in this case. On the other hand, if more than $n$ vehicles have been explored, the next target in the sorted list based on savings for vehicle $i$ is considered.
    \end{enumerate}
\end{enumerate}
Similar to the previous neighborhood, the search for a better solution in this neighborhood is continued till a better solution is obtained or all targets from the maximal vehicle have been explored for a swap. The sequence of steps in this neighborhood is depicted in Fig.~\ref{fig: swap_neighborhood}.

\subsubsection{Neighborhood 3: multi-target swap} \label{subsubsect: multi-target}

\begin{figure}[htb!]
    \centering
    \subfigure[Vehicle tours before removing targets $t_1$ and $t_2$]{\includegraphics[width = 0.32\textwidth]{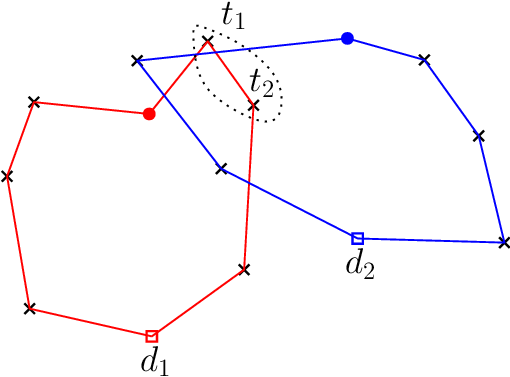}} \hfill
    \subfigure[Vehicle tours after inserting $t_1$ and $t_2$ into ``blue" vehicle]{\includegraphics[width = 0.32\textwidth]{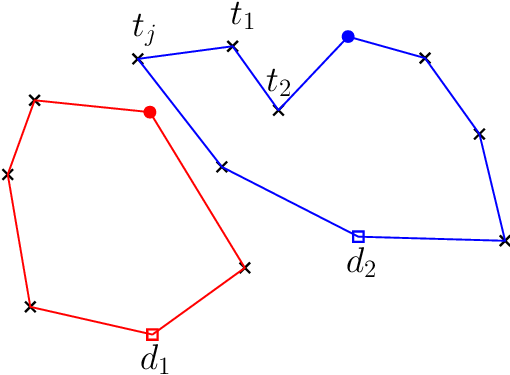}} \hfill
    \subfigure[Vehicle tours after removing target $t_j$ from ``blue" vehicle and inserting in ``red" vehicle]{\includegraphics[width = 0.32\textwidth]{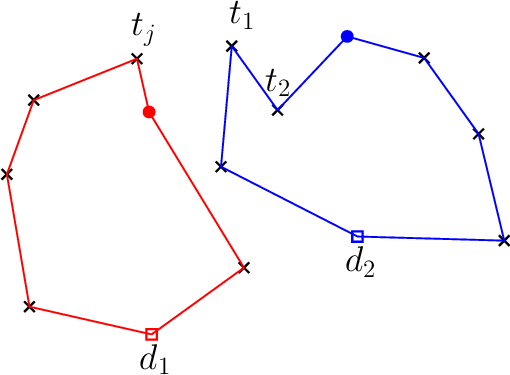}}
    \caption{Illustration of multi-target swap neighborhood for exchange from maximal vehicle (in red) and another vehicle (in blue). Targets colored in black are common targets, and shaded discs represent targets with vehicle-target assignments.}
    \label{fig: neighborhood_3}
\end{figure}

In the multi-target swap neighborhood, a set of targets from the maximal vehicle $i$ is attempted to be swapped with a set of targets from another vehicle, as shown in Fig.~\ref{fig: neighborhood_3}. This study considers two types of multi-target swap neighborhoods: fixed-structure and variable-structure. In the fixed-structure multi-target swap neighborhood, a set of $m$ targets from the maximal vehicle is attempted to be swapped with a set of $m - 1$ or $m$ targets\footnotemark\, from another vehicle. However, since a heterogeneous fleet of vehicles is considered, exchanging a similar number of targets might not substantially improve the solution since one vehicle might have a higher speed than another. Therefore, a variable-structure multi-target swap neighborhood is also considered in this study, wherein multi-target swaps of varying sizes are considered.

\footnotetext{In this study, exchange with $m - 1$ targets was also considered to account for particularly swapping two targets from the maximal vehicle with one target in another vehicle, since to attain such a solution, a mix of switch and swap neighborhoods will have to be utilized. However, for cases wherein a target switch or a target swap leads to a solution with a higher objective cost, such a solution cannot explored.}

\textbf{Remark:} A figure for the outline for the implementation of the multi-target swap neighborhood is not provided for brevity due to its similarity to the implementation of the target swap neighborhood given in Fig.~\ref{fig: swap_neighborhood}.

\paragraph{Fixed-structure multi-target swap:}

In the fixed-structure multi-target swap neighborhood, it is desired to exchange a set of $m$ targets from the maximal vehicle with $m - 1$ or $m$ targets from another vehicle. To this end, targets that are pre-assigned to vehicle $i,$ i.e., in the set $R_i,$ are first excluded. Noting that the number of combinations of targets that can be removed can be large, especially with increasing $m$, the neighborhood is restricted to consider a set of $m$ targets that are subsequent targets in the tour, or only separated by targets in $R_i$ in the tour. The savings associated with removing a group of $m$ targets from the maximal vehicle is then computed.

It should be noted that the savings associated with removing a group of targets is given by the sum of savings of removing targets one after the other from the considered vehicle's tour. For example, consider removing targets $u_l^i, u_{l + 1}^i$ from the maximal vehicle's tour given by $d_i, u_1^i, u_2^i, \cdots, u_{l - 1}^i, u_l^i, \cdots, u_{l + m - 1}^i, u_{l + m}^i, \cdots, d_i.$ First, $u_l^i$ is removed and its savings is computed using Eq.~\eqref{eq: definition_savings_metric}, which corresponds to removing edges $(u_{l - 1}^i, u_l^i)$ and $(u_l^i, u_{l + 1}^i)$ and adding edge $(u_{l - 1}^i, u_{l + 1}^i)$ to construct a feasible tour. Then, the savings associated with $u_{l + 1}^i$ is computed using Eq.~\eqref{eq: definition_savings_metric} using the current tour, wherein $u_{l - 1}^i$ is the target preceding it in the tour. Therefore, the computed savings corresponds to removing edges $(u_{l - 1}^i, u_{l + 1}^i)$ and $(u_{l + 1}^i, u_{l + 2}^i)$ and adding the edge $(u_{l - 1}^i, u_{l + 2}^i).$ It can be verified that such a computation is independent of the order of removal of targets from the tour. 

\begin{figure}[htb!]
    \centering
    \subfigure[Portion of tour before $2-$opt swap of edges \label{subfig: before_2_opt_swap}]{\includegraphics[width = 0.45\textwidth]{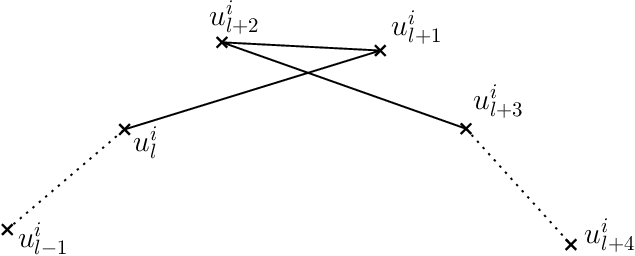}} \hfill
    \subfigure[Portion of tour after $2-$opt swap of edges \label{subfig: after_2_opt_swap}]{\includegraphics[width = 0.45\textwidth]{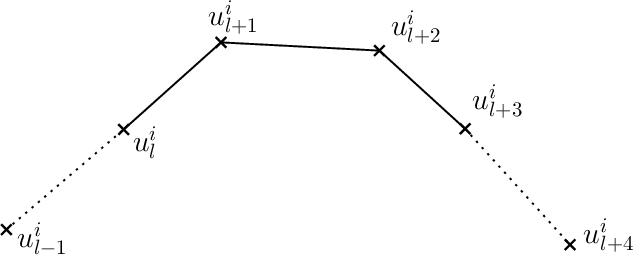}}
    \caption{Portion of tour before and after $2-$opt swap of edges} \label{fig: 2-opt_swap}
\end{figure}

Using the computed savings for each set of $m$ targets from the maximal vehicle, the target sets are sorted in the decreasing order of savings. Then, the following set of steps are performed to improve the feasible solution:
\begin{enumerate}
    \item First, it is desired to pick a vehicle in which the set of targets must be inserted. Such a vehicle is selected based on the insertion cost since this metric was observed to provide the best improvement in solution in a computationally efficient manner for the target switch and swap neighborhoods. (The results for the same will be presented in Section~\ref{sect: results}.) It is now desired to define the insertion cost for inserting a set of targets from the maximal vehicle, say $u_l^i, u_{l + 1}^i, \cdots, u_{l + m - 1}^i$, in another vehicle. In the fixed-structure multi-target swap neighborhood, the considered targets for insertion are inserted as a group/set between two targets in another vehicle to reduce the number of choices present.

    Suppose the insertion cost for vehicle $j$ is desired to be computed. Let vehicle $j$'s tour be given by $d_j, u_1^j, \cdots, u_p^j, u_{p + 1}^j, \cdots, d_j$, and the insertion cost for inserting the considered group of targets between $u_p^j$ and $u_{p + 1}^j$ is desired to be computed. First, if the number of targets in vehicle $j$'s tour that are not preassigned to vehicle $j$ is less than $m - 1,$ such a vehicle is not considered since the vehicle does not have sufficient number of targets to perform a swap with vehicle $i$. Noting that the group of targets can be inserted in two orientations, wherein either $u_p^j$ is connected to the first target in the group ($u_l^i$) and $u_{p + 1}^j$ is connected to the last target in the group ($u_{l + n - 1}^i$) in the first orientation, whereas it is vice-versa in the second orientation, the insertion cost is given by
    \begin{align} \label{eq: insertion_cost_fixed_structure}
    \begin{split}
        &\text{insertion cost}_{(u_l^i, u_{l + 1}^i, \cdots, u_{l + m - 1}^i), j, (u_p^j, u_{p + 1}^j)} \\
        &= \min\Bigg(\frac{dist (u_p^j, u_l^i) + \sum_{q = 1}^{m - 1} dist (u_{l + q - 1}^i, u_{l + q}^i) + dist (u_{l + n - 1}^i, u_{p + 1}^j)}{v_j}, \\
        & \qquad\quad\quad \frac{dist (u_p^j, u_{l + m - 1}^i) + \sum_{q = 1}^{n - 1} dist (u_{l + q - 1}^i, u_{l + q}^i) + dist (u_{l}^i, u_{p + 1}^j)}{v_j} \Bigg) \\
        & \quad\, - \frac{dist (u_p^j, u_{p + 1}^j)}{v_j}.
    \end{split}
    \end{align}
    Since the set of targets can be inserted in any location in vehicle $j$'s tour, the insertion cost for the considered set of targets for vehicle $j$ corresponds to the location that yields the minimum cost. In this manner, the vehicle that yields the minimum cost of insertion is selected, and the set of targets is inserted in the considered vehicle at the location with an orientation that yields the least insertion cost.
    \item Suppose vehicle $j$ is the vehicle that yields the minimum cost of insertion, whose tour was updated with the set targets removed from vehicle $i$. A $2-$opt swap of edges is performed to improve the tour of vehicles $i$ and $j$ in a computationally efficient manner \cite{2_opt, dubins_tsp_heuristic_nayak}. A depiction of a portion of the vehicle tour before and after a $2-$opt swap is shown in Fig.~\ref{fig: 2-opt_swap}. 
    \item Now, it is desired to pick a set of $m$ targets or $m - 1$ targets from vehicle $j$ and insert them in vehicle $i$. To this end, targets that are pre-assigned to vehicle $j,$ i.e., in $R_j,$ and the set of targets inserted in vehicle $j$ from vehicle $i$ are first excluded. Furthermore, sequential target sets in the vehicle's tour are considered to reduce the number of solutions explored, and the savings obtained by removing from vehicle $j$ and the insertion cost for inserting in vehicle $i$ are computed using the same method explained beforehand. Now, it is desired to sort the target sets from vehicle $j$ to determine if a better solution can be obtained. In this regard, two sorting strategies are employed:
    \begin{enumerate}
        \item Target sets in vehicle $j$ are sorted in the increasing order of the insertion cost associated with inserting in vehicle $i$, since target sets that are ``cheap'' for vehicle $i,$ the maximal vehicle, to cover is desired.
        \item Target sets in vehicle $j$ are sorted in the decreasing order of the difference between the savings obtained from vehicle $j$ minus the insertion cost associated with inserting in vehicle $i$. Such a sorting strategy is considered since it is desired to consider target sets that are ``expensive'' for vehicle $j$ to cover and are ``cheap'' for vehicle $i$ to cover.
    \end{enumerate}
    \item Having sorted the target sets in vehicle $j$ through one of the strategies, each candidate target set from the sorted list is attempted to be removed from vehicle $j$ and inserted in vehicle $i$. In this regard, the location for insertion in vehicle $i$ and orientation are chosen based on the least insertion cost, defined in Eq.~\eqref{eq: insertion_cost_fixed_structure}, and a $2-$opt swap of edges is performed for both vehicle $i$ and vehicle $j$. Using the tour costs of vehicles $i$ and $j,$ which are kept track of through the savings and insertion costs (similar to Neighborhood~2) and a $2-$opt swap of edges, the maximum tour cost among the vehicles is compared with the incumbent tour cost. If a better solution is obtained, LKH will be used to optimize the tours of both vehicles, and the incumbent solution will be updated. If not, the next target set in the sorted list of targets from vehicle $j$ is explored. Since the number of target sets can be large, the number of candidate target sets explored in vehicle $j$ is restricted based on a parameter $n$. If the number of target sets explored in vehicle $j$ exceeded the maximum number of candidate target sets or all target sets in vehicle $j$ have been covered, the next target set for removal from vehicle $i,$ the maximal vehicle, is considered.
\end{enumerate}
The same sequence of steps is continued till a better solution is obtained, or all target sets from vehicle $i$ have been considered.

\textbf{Remark:} Unlike the switch and the swap neighborhoods, multiple vehicles are not considered for inserting a selected set of targets from the maximal vehicle since the multi-target swap neighborhood is a much larger neighborhood. In this regard, to reduce the number of candidate solutions explored and, therefore, the computation time, only one vehicle is considered corresponding to a selected set of targets from the maximal vehicle. The same argument applies to the variable-structure multi-target swap neighborhood, the details of which follow.

\paragraph{Variable-structure multi-target swap:}

\begin{figure}[htb!]
    \centering
    \subfigure[Target cluster with high removal ratio \label{subfig: variable_structure_1}]{\includegraphics[width = 0.45\textwidth]{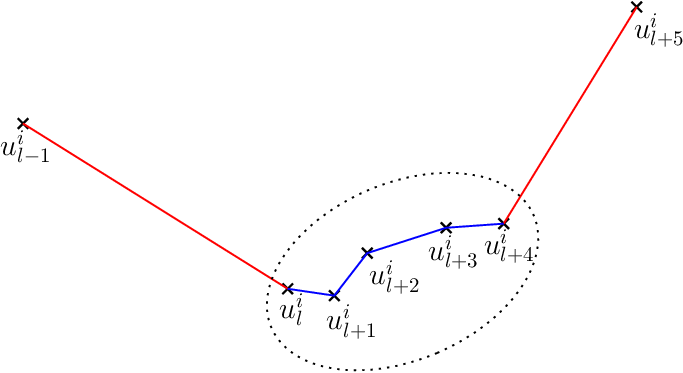}} \hfill
    \subfigure[Target cluster with lower removal ratio \label{subfig: variable_structure_2}]{\includegraphics[width = 0.45\textwidth]{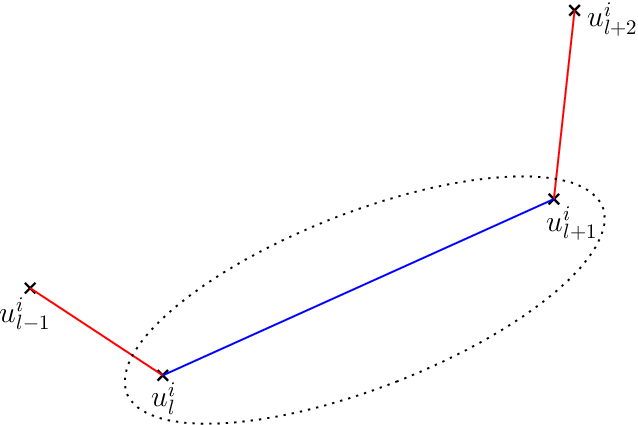}}
    \caption{Depiction of removal ratio for variable-structure multi-target swap neighborhood} \label{fig: variable_structure_multi_target_swap}
\end{figure}

In the variable-structure multi-target swap, exchanging a set of at most $m$ targets from the maximal vehicle with a target set of at most $m$ targets with another vehicle is desired. In this regard, target sets of varying sizes from $2$ to $m$ are considered from the maximal vehicle. Similar to the fixed-structure multi-target swap variant, targets that are pre-assigned to vehicle $i$ (in set $R_i$) are removed, and target sets that are sequential in the vehicle tour are considered (after removing the pre-assigned targets). The motivation behind this neighborhood is to account for the heterogeneity of vehicles in the neighborhood. However, unlike the fixed-structure multi-target swap neighborhood, the savings metric alone cannot be used since target sets of size $2$ will predominantly be favored over target sets of a larger size. In this regard, the sorting of targets is performed using a ``removal ratio''. For vehicle $i$, whose tour is given by $d_i, u_1^i, \cdots, u_l^i, \cdots, u_{l + p - 1}^i, \cdots, d_i$, the ``removal ratio" associated with removing a target set $u_l^i, \cdots, u_{l + p - 1}^i$ in vehicle $i$'s tour is defined as
\begin{align*}
    \text{removal ratio}_{u_l^i, \cdots, u_{l + p - 1}^i} = \frac{dist (u_{l - 1}^i, u_l^i) + dist (u_{l + p - 1}^i, u_{l + p}^i)}{\sum_{q = 1}^{p - 1} dist (u_{l + q - 1}^i, u_{l + q}^i)}.
\end{align*}
The defined removal ratio can be used to identify target clusters that are far from the rest of the tour; hence, a target cluster depicted in Fig.~\ref{subfig: variable_structure_1} will have a higher removal ratio compared to a target cluster depicted in Fig.~\ref{subfig: variable_structure_2}. The target sets are sorted in the decreasing order of the removal ratio since it is desired to remove target clusters that are expensive for vehicle $i$ to visit and are explored one at a time for removal and insertion into another vehicle.

The following set of steps are then repeated to determine if a better solution can be obtained:
\begin{enumerate}
    \item Consider inserting the target set $u_l^i, \cdots, u_{l + p - 1}^i$, which is the higher removal ratio in the maximal vehicle $i,$ where $p \leq m$. The vehicle chosen for inserting the considered target set is based on the insertion cost. It must be recalled that the insertion cost of inserting a set of targets was defined for the fixed-structure multi-target swap neighborhood. The expression for inserting in vehicle $j$ between targets $u_p^j$ and $u_{p + 1}^j$ is given in Eq.~\eqref{eq: insertion_cost_fixed_structure}. While such an insertion is computationally efficient, one other definition of insertion cost is explored, which is considered a variant in this neighborhood. The alternate variant considers inserting the targets from the target set one after the other in another vehicle's tour (denoted by $j$). The steps for the same are as follows:
    \begin{enumerate}
        \item Consider the first target in the target set to be inserted ($u_l^i$) and insert the target in the location with the least cost of insertion in vehicle $j$'s tour, where the cost of inserting between two targets in vehicle $j$ is defined in Eq.~\eqref{eq: insertion_cost_single_vehicle_definition}. Initialize the insertion cost for the set to be zero when the first target in the target set ($u_l^i$) is considered.
        \item Remove the considered target from the target set and update the insertion cost of the target set with the previously computed insertion cost for the target.
        \item Repeat the above two steps till all targets from the target set have been considered.
    \end{enumerate}
    \item A $2$-opt swap of edges is performed for the maximal vehicle and the vehicle in which the target set is inserted (denoted as vehicle $j$). It is now desired to remove a set of targets from vehicle $j$. To this end, noting that target sets of differing sizes need to be considered, and noting that it is desired to remove sets that are ``expensive" for vehicle $j$ and ``cheap" for vehicle $i,$ the target sets in vehicle $j$ of size at most $m$ are sorted in the decreasing order of the ratio of savings associated with removing from vehicle $j$ and inserting in vehicle $i.$ It should be noted that in the first variant, insertion cost based on inserting as a group is computed (similar to the fixed-structure multi-target swap neighborhood), whereas, in the second variant, insertion cost based on recursive insertion is considered.
    \item Each target set in the sorted list from vehicle $j$ is considered for removal and insertion into vehicle $i.$ A $2-$opt swap of edges is then performed, and based on estimated tour costs, it is determined if a better solution is obtained or not. If a better solution is obtained, then LKH is used to optimize the tours. If not, then the next target set in vehicle $j$ is explored if the maximum number of candidate sets $n$ has not been explored (similar to the fixed-structure neighborhood). If the maximum number of sets have been explored, then the next target set in vehicle $i$ is explored.
\end{enumerate}
The same sequence of steps is performed till a better solution is obtained or all target sets from vehicle $i$ have been considered.

\subsection{Perturbation}

\begin{figure}
    \centering
    \includegraphics[width=0.3\linewidth]{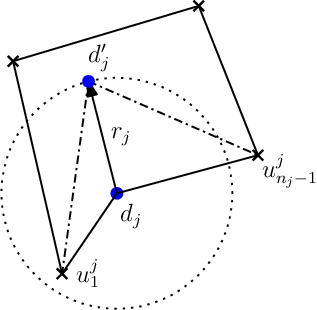}
    \caption{Depiction of perturbation of the depot of the $j\textsuperscript{th}$ vehicle \cite{AIAA_modified_MD}}
    \label{fig: perturbation_j_vehicle}
\end{figure}

Through neighborhood searches, a local minimum would be obtained. The perturbation step is used to break from this local minimum, wherein the solution is perturbed at random to explore a potentially new solution. To this end, the perturbation step adopted in \cite{AIAA_modified_MD} is utilized.
In the first stage of the perturbation step, the depot of the $j\textsuperscript{th}$ vehicle ($d_j$) is perturbed by a distance $r_j$ along a random angle $\theta_j.$ Here, the distance $r_j$ is selected to be 
\begin{align*}
    r_j = \frac{1}{2} \left(\frac{dist (u_{n_j - 1}^j, d_j) + dist (d_j, u_1^j)}{v_j} \right),
\end{align*}
which is the average cost of the edges that connect $d_j$ to the first target ($u_1^j$) and the last target ($u_{n_j - 1}^j$) visited by vehicle $j$ in its tour in the incumbent solution. A depiction of such a perturbation is shown in Fig.~\ref{fig: perturbation_j_vehicle}.

Followed by such a perturbation for all vehicles (i.e., $j = 1, 2, \cdots, k$), a new graph is obtained with new locations for all depots. Using the same allocation of targets to vehicles from the incumbent solution, the tours in the new graph are constructed using LKH. A local search is then performed using the local search method explained before to obtain a local minimum in this new graph. Using the allocated targets for each vehicle in this new graph, LKH is utilized on the original graph (with the initial actual depot locations) to obtain new potential tours for each vehicle. If an improved solution is obtained compared to the incumbent, then the incumbent solution is updated, and the perturbation step is stopped. If not, the perturbation step is repeated, wherein in the following perturbation, the angle of perturbation of the $j\textsuperscript{th}$ depot is $144^\circ$ more than $\theta_j.$ At most five perturbations are performed with no improvement in each perturbation step before it is stopped (refer to Fig.~\ref{fig: overview_heuristic}).

\section{Results} \label{sect: results}

In this section, computational results are presented to determine the best among the many variations proposed in the heuristic in this study. The computations presented in this section are based on implementations in Python~3.9 on a laptop with AMD Ryzen $9$ $5900$HS CPU running at $3.30$ GHz with $16$ GB RAM. Furthermore, the linear program for the LP-based load balancing method for generating the initial feasible solution was solved using Gurobi \cite{gurobi} in Python. The organization of the presented results in this section is as follows:
\begin{enumerate}
    \item The developed set of 43 instances in \cite{MD_algorithm} for the homogeneous problem is first utilized to generate a set of 128 instances for the heterogeneous problem.
    \item The two methods considered for generating an initial feasible solution, presented in Section~\ref{subsect: initialisation_step}, are compared over the generated set of heterogeneous instances.
    \item The various variations in the switch and swap neighborhood are then analyzed over the set of 128 instances, along with the two different methods of generating the initial feasible solution. The initial feasible solution methods were also considered in this analysis since, from the comparison of the objective value of the initial solutions generated from the two methods, the better-performing method could not be identified. Through extensive computational results and comparisons, one among the three metrics for the target switch and swap neighborhoods, the best-performing value of $n$, wherein $n$ determines the number of vehicles considered for the target switch/swap, and the best method for constructing the initial feasible solution are picked.
    \item Having identified the best switch and swap neighborhood, a list of the most sensitive instances is identified to perform the multi-target neighborhood analysis. Through an identified set of 36 instances out of 128 instances, the variations in the multi-target swap are analyzed, and the best-performing metric and parametric value for the multi-target swap neighborhood in terms of the objective value and computation time are identified.
    \item Finally, a summary of results is provided, wherein the heuristic with the best switch and swap neighborhood and the heuristic with the best switch and swap neighborhood along with the best multi-target swap neighborhood are benchmarked against the heuristic in \cite{AIAA_modified_MD}, termed the modified MD algorithm. Furthermore, the tours obtained from the benchmarking algorithm and the two heuristic variants that are finally considered will be shown for a sample instance.
\end{enumerate}

\textbf{Remark:} From the literature surveyed before, it can be observed that only one heuristic has been proposed to address the considered heterogeneous min-max problem \cite{AIAA_modified_MD}. Hence, the presented heuristic in this paper and its variations will be benchmarked against this algorithm.
    


\subsection{Generation of Heterogeneous Instances}

In \cite{MD_algorithm}, a total of 43 instances were generated with a number of targets ranging from 10 to 500 and a number of vehicles ranging from 3 to 20. Furthermore, the location of each target was generated based on different probability distributions in \cite{MD_algorithm} for different instances. Due to the wide variety of datasets considered for the homogeneous variant of the problem, it is desired to utilize the target and depot locations in each instance for the heterogeneous variant of the problem. In these instances, structural heterogeneity is introduced by randomly assigning a speed in the set $\{1, 1.25, 1.5, 1.75, 2\}$ to each vehicle. Additionally, three cases of functional heterogeneity were considered for each such instance:
\begin{itemize}
    \item No functional heterogeneity. In this case, all targets in an instance can be covered by any vehicle, i.e., $R_i = \emptyset \, \forall \, i = 1, 2, \cdots, k.$
    \item Functional heterogeneity with three targets assigned per vehicle. In this case, for the first vehicle, the available list of targets that have not been assigned is first sorted in terms of the distance from the depot of the vehicle, and the closest target, median distance target, and farthest target are assigned to the vehicle. The assigned targets are then removed from the common list of targets that can be covered by any vehicle, and the same process is repeated for vehicles $2, 3, \cdots, k.$
    \item Functional heterogeneity with five targets assigned per vehicle. Similar to the previous case, the available targets are first sorted based on the distance from the depot for the first vehicle. The two closest targets, two farthest targets, and the median distance target are assigned to the vehicle and removed from the common list of targets that can be covered by any vehicle. The same process is repeated for vehicles $2, 3, \cdots, k.$
\end{itemize}
In this regard, 43 instances were obtained for each of the zero-target allocation and three-target allocation cases using the 43 instances in \cite{MD_algorithm}. For the case with five-target allocations, 42 instances were obtained since, for the first instance in \cite{MD_algorithm}, the ratio of targets to vehicles was less than five. Therefore, a total of 128 instances are considered for the heterogeneous problem.

\subsection{Comparison of Initial Feasible solution}

\begin{figure}[htb!]
    \centering
    \subfigure[Comparison of objective value for the case with no vehicle-target assignment \label{subfig: comparison_no_vehicle_target_assignment}]{\includegraphics[width = 0.45\linewidth]{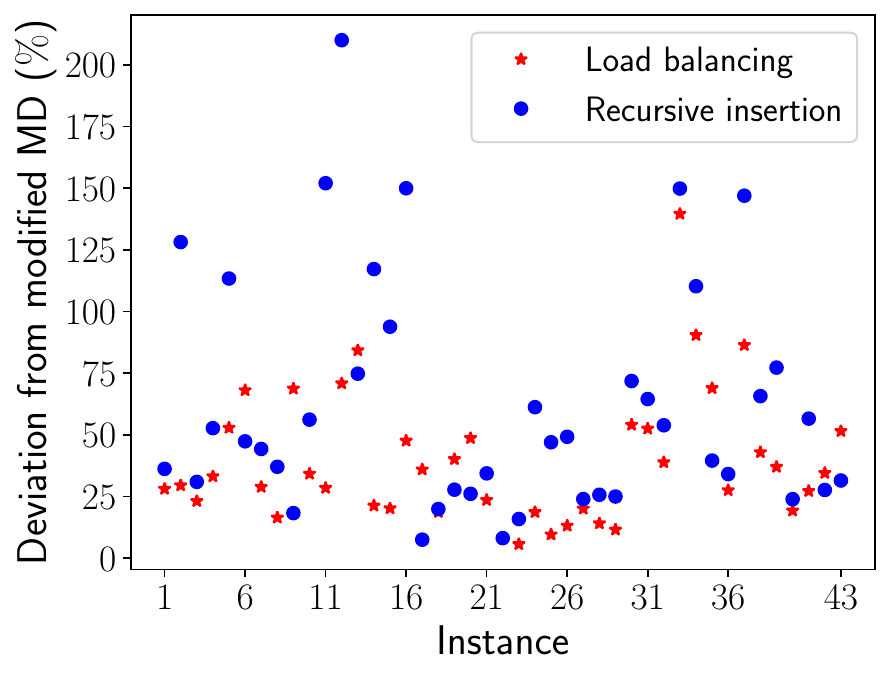}}\hfill
    \subfigure[Comparison of objective value for the case with three targets assigned per vehicle \label{subfig: comparison_three_vehicle_target_assignment}]{\includegraphics[width = 0.45\linewidth]{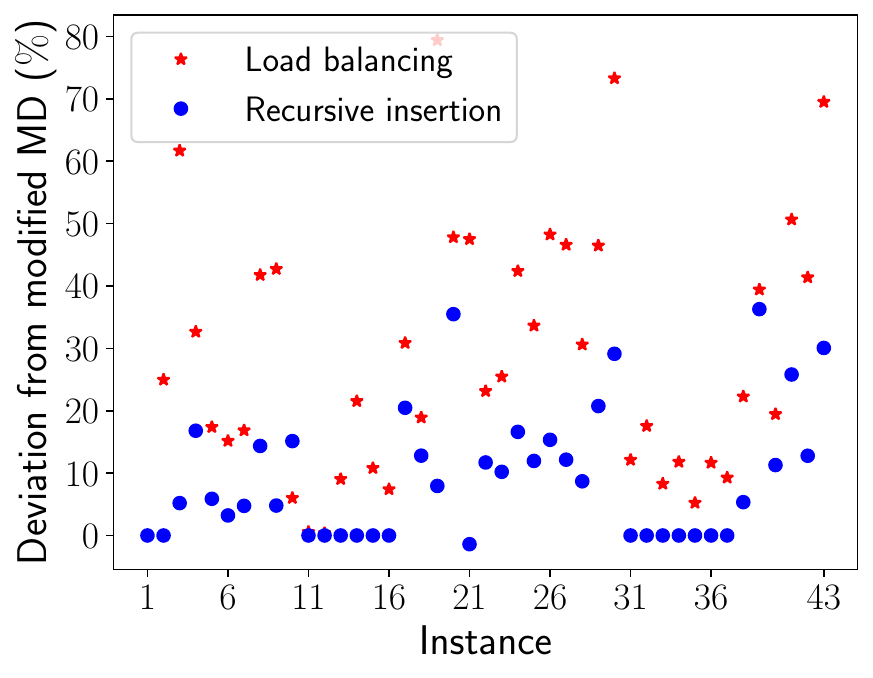}}\hfill
    \subfigure[Comparison of computation time for the case with no vehicle-target assignment \label{subfig: comparison_no_vehicle_target_assignment_comp_time}]{\includegraphics[width = 0.45\linewidth]{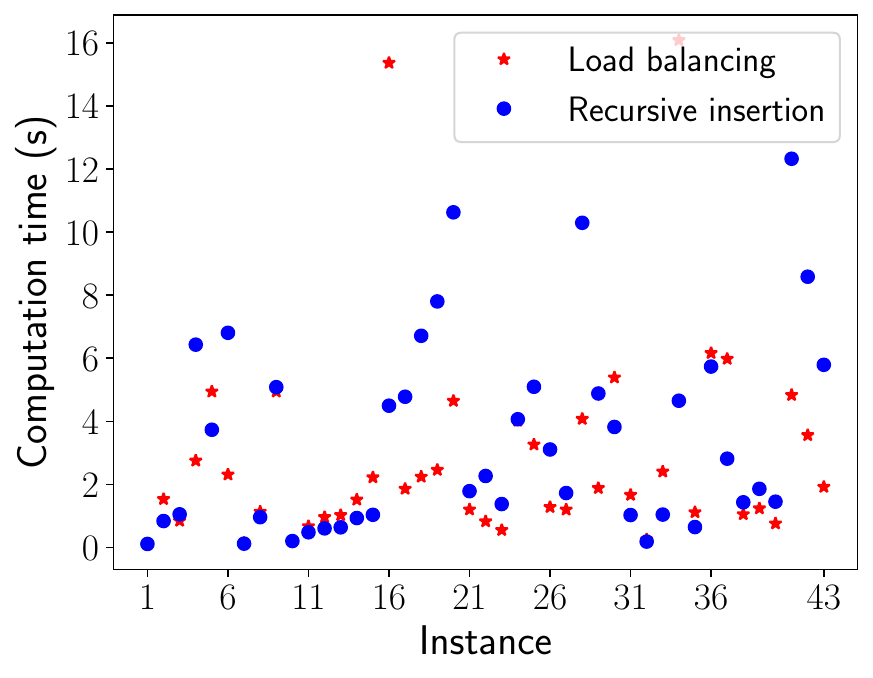}} \hfill
    \subfigure[Comparison of computation time for the case with three targets assigned per vehicle \label{subfig: comparison_three_vehicle_target_assignment_comp_time}]{\includegraphics[width = 0.45\linewidth]{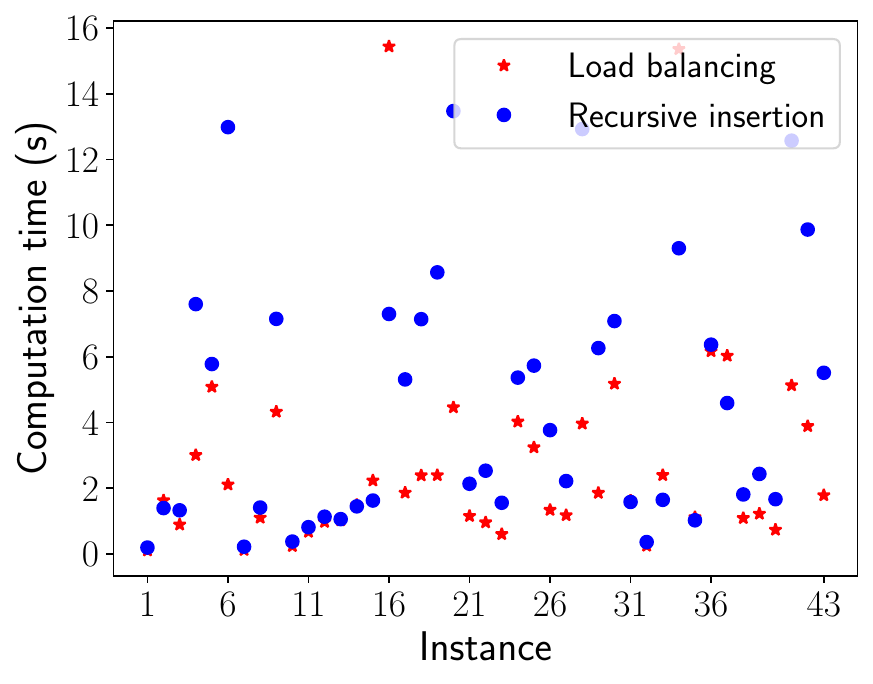}}
    \caption{Comparison of deviation of initial feasible solution from benchmarking algorithm \cite{AIAA_modified_MD} and computation time} \label{fig: initial_feasible_soln_comparison}
\end{figure}

Using the generated heterogeneous instances, it is first desired to determine the best method among the two methods discussed in Section~\ref{subsect: initialisation_step} to generate an initial feasible solution. The comparison of the objective value of the initial feasible solution generated for the instances with zero-target allocation and three-target allocation is shown in Figs.~\ref{subfig: comparison_no_vehicle_target_assignment} and \ref{subfig: comparison_three_vehicle_target_assignment}, respectively. Furthermore, the comparison of the computation time to generate the initial feasible solution and utilize LKH to optimize each vehicle's tour is shown for the zero-target and three-target cases in Figs.~\ref{subfig: comparison_no_vehicle_target_assignment_comp_time} and \ref{subfig: comparison_three_vehicle_target_assignment_comp_time}, respectively. It should be noted that the comparison for the five-target case is not shown due to its similarity to the three-target case. A summary of the results is also provided in Table~\ref{tab: comparison_initial_soln}. From these figures and table, it can be observed that
\begin{itemize}
    \item The load balancing method yields a better initial solution for the zero-target cases, whereas the recursive insertion method outperforms the load balancing method for the three-target and five-target cases.
    \item Both methods take a similar computation time to yield the initial solution, with the load balancing method yielding the initial solution marginally faster than the recursive insertion method by a couple of seconds.
\end{itemize}
Due to the mixed results obtained with regard to the objective value and the negligible differences in the computation time, the best choice of method for generating the initial solution is not evident. Due to this reason, the two methods are considered variations in the heuristic and are compared in terms of the solution obtained from the heuristic with a switch and swap neighborhood in Section~\ref{subsubsect: switch_swap_analysis}.

\begin{table}[htb!]
    \centering
    \caption{Comparison of the number of instances in which the initial feasible solution generated from load balancing (LB) method was better, equal, or worse than the solution from the recursive insertion method}
    \label{tab: comparison_initial_soln}
    \begin{tabular}{ccccccc} \toprule
    \textbf{Functional} & \multicolumn{3}{c}{\textbf{Soln. comparison}} & \multicolumn{3}{c}{\textbf{Computation time comparison}} \\ \cmidrule{2-4} \cmidrule{5-7}
    \textbf{heterogeneity} & LB better & Equal & LB worse & LB better & Equal & LB worse \\ \midrule
    No target & 34 & 0 & 9 & 24 & 0 & 19 \\
    Three targets & 1 & 1 & 41 & 34 & 0 & 9 \\
    Five targets & 1 & 1 & 40 & 32 & 0 & 10 \\ \bottomrule
    \end{tabular}
\end{table}

\subsection{Target Switch and Swap Neighborhood Analysis} \label{subsubsect: switch_swap_analysis}

It must be recalled that a total of three metrics were considered in the definition of the switch neighborhood and swap neighborhood, which are Neighborhoods $1$ and $2$, respectively. Furthermore, a parameter $n$ was considered, which denotes the number of vehicles considered for switching or swapping a target from the maximal vehicle. In this regard, it is desired to identify the best metric choice and $n$ value that provides the best trade-off between objective value improvement and computation time. Furthermore, it is also desired to identify the method of initial solution construction that yields the best objective value. Hence, two variations in the construction method, three variations in the metric, and $n \in \{1, 2, 3\}$ are considered.

\begin{figure}[htb!]
    \centering
    \subfigure[Variation due to metric on objective value \label{subfig: variation_metric_obj_val}]{\includegraphics[width = 0.48\textwidth]{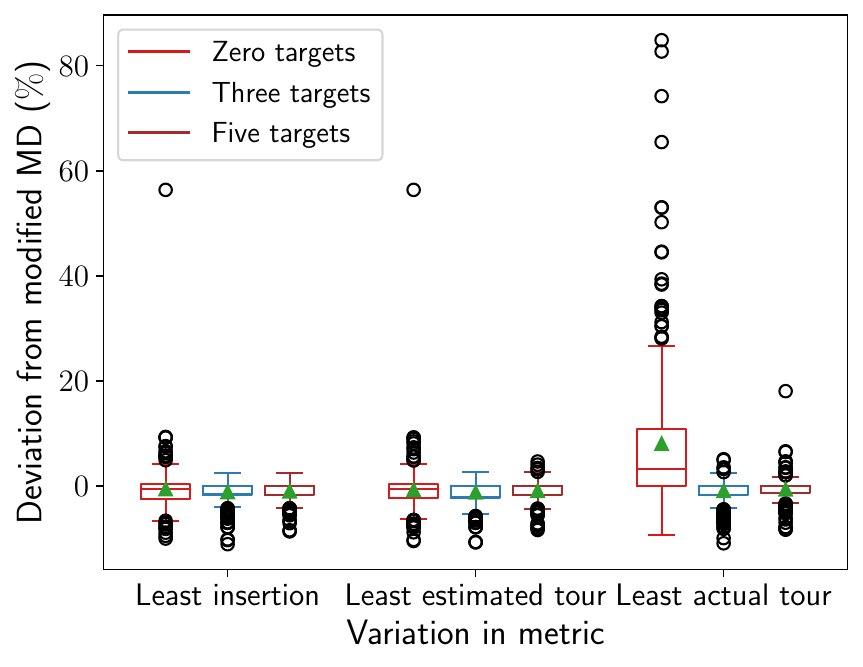}} \hfill
    \subfigure[Variation in objective value based on construction method for zero target allocation instances \label{subfig: variation_objective_val_zero_targ}]{\includegraphics[width = 0.48\textwidth]{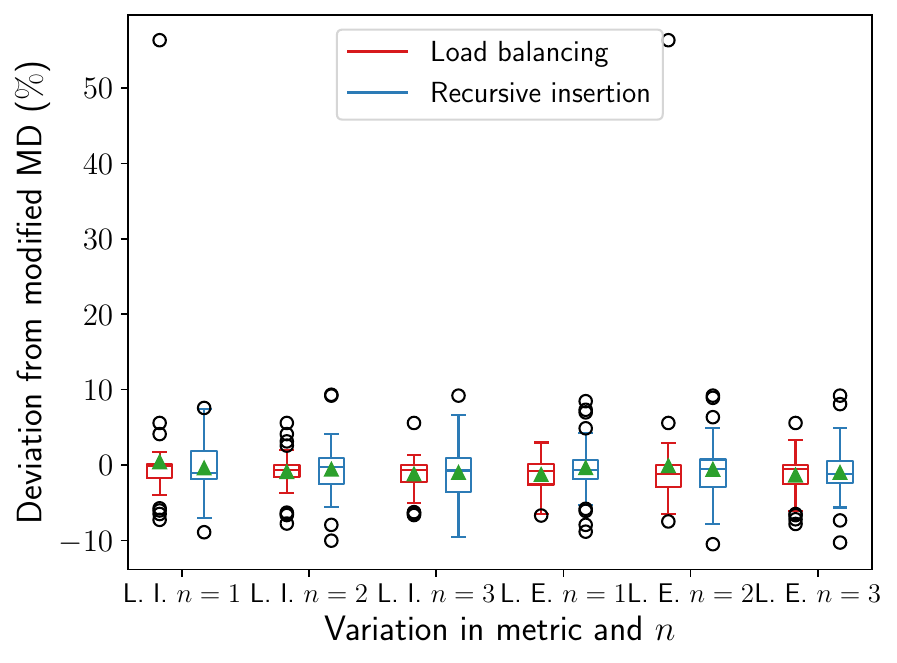}}
    \subfigure[Variation in objective value based on construction method for three target allocation instances \label{subfig: variation_objective_val_three_targ}]{\includegraphics[width = 0.48\textwidth]{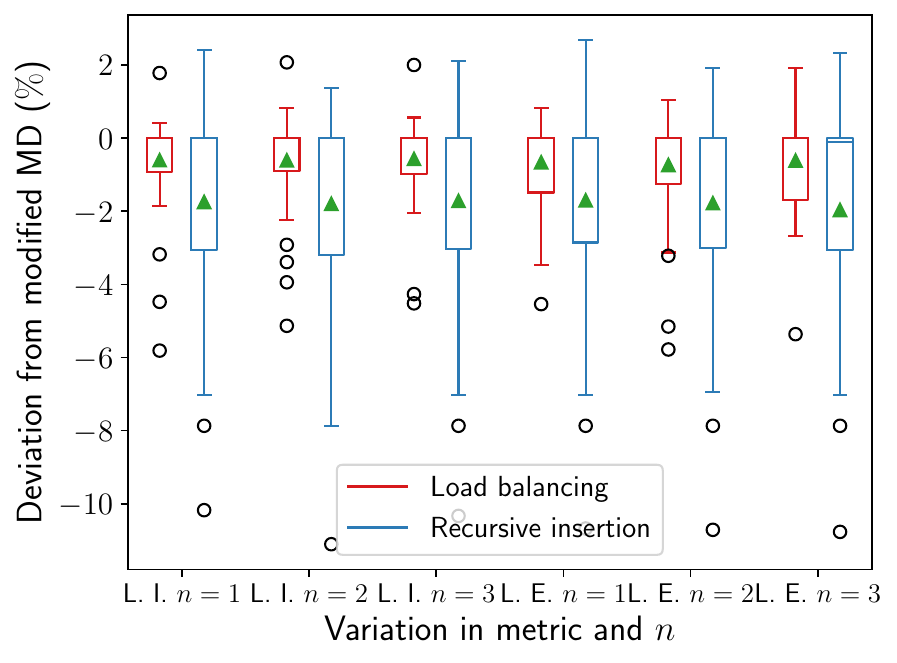}}
    \caption{Analysis of switch and swap neighborhood} \label{fig: analysis_switch_swap_neighborhood}
\end{figure}

\begin{figure}[htb!]
    \centering
    \includegraphics[width=\linewidth]{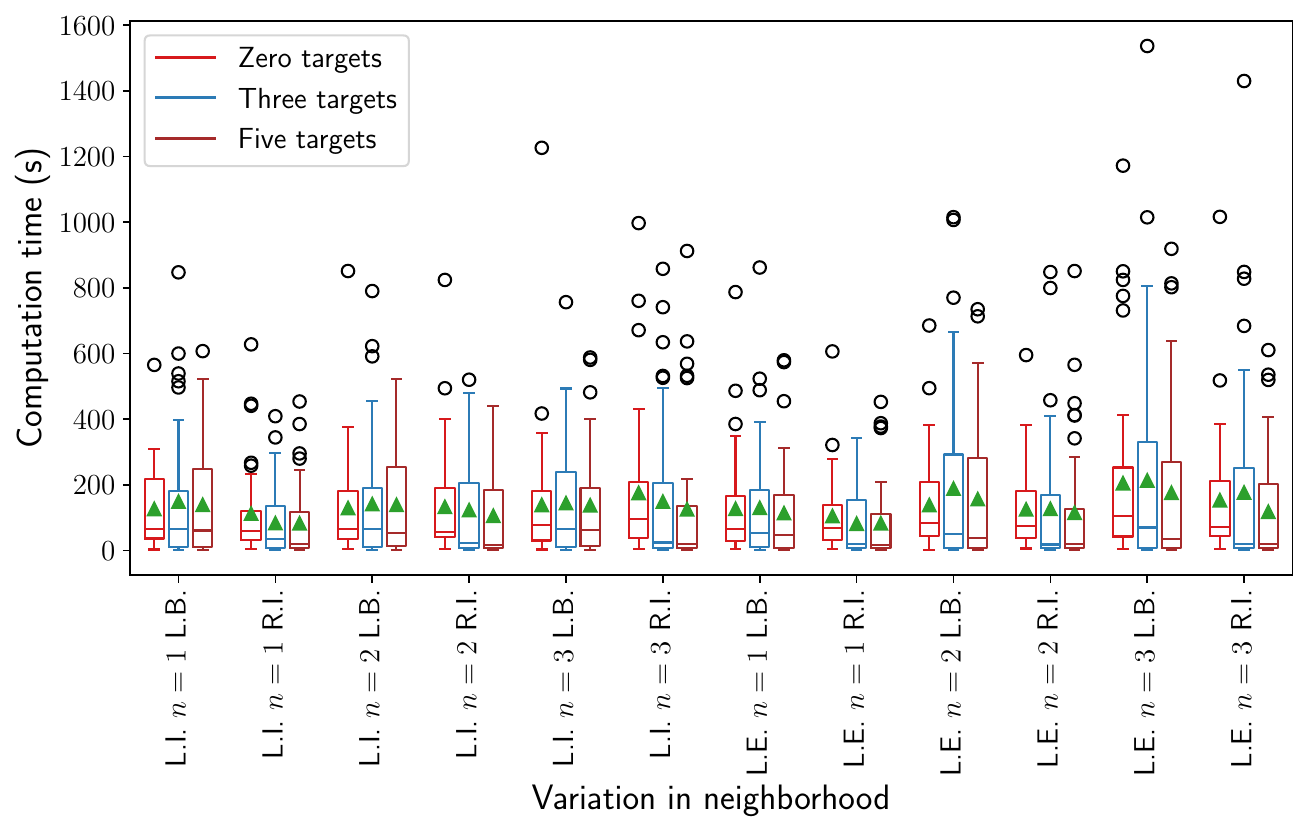}
    \caption{Box plot depicting the computation time for the two chosen metrics (least insertion, abbreviated as L.I., and least estimated tour, abbreviated as L.E.) with three different $n$ values and two construction methods, which are load balancing (abbreviated as L.B.) and recursive insertion (abbreviated as R.I.) for the switch and swap neighborhoods for the three types of instances (along with the mean time depicted using a triangle)}
    \label{fig: computation_time_box_plot_switch_swap}
\end{figure}

The variations were evaluated over the $128$ instances generated for the heterogeneous problem.
To account for randomness, which primarily arises from the perturbation step, the heuristic was run three times for each instance, and the best solution obtained, along with the corresponding computation time, was recorded. The observations 
obtained are as follows:
\begin{enumerate}
    \item The percentage deviation of the objective values from the modified MD algorithm for the three different metrics are depicted in a box diagram in Fig.~\ref{subfig: variation_metric_obj_val}, wherein the mean deviation is also depicted (using a triangle). From this figure, it can be observed that the three metrics produce similar results for the cases of three-target and five-target allocations. However, the least actual tour cost metric produces poor results in the case of zero-target allocation. In this case, the median difference from the modified MD algorithm is positive; furthermore, from the results, it can be observed that the maximum deviation is more than $20\%$ for many instances. Hence, the least actual tour cost metric for the neighborhood will not be considered.
    \item For each of the other two metrics considered, the percentage deviation from the modified MD algorithm for each value of $n$ is depicted for the zero-target instances and three-target instances\footnotemark\, in Figs.~\ref{subfig: variation_objective_val_zero_targ} and \ref{subfig: variation_objective_val_three_targ}, respectively.
    It can be observed that the recursive insertion method for constructing the initial solution yields a better solution for the heuristic than the load balancing method. In particular,
    \footnotetext{It must be noted here that the deviation plot for the five-target case is very similar to the three-target case and is not shown for brevity.}
    \begin{enumerate}
        \item For the zero-target allocation instances, the solution obtained from the heuristic using load balancing yields outliers for one of the $n$ values for both metrics. Additionally, the minimum deviation obtained from the recursive insertion method is better than the load balancing method by about $2$ to $3\%$.
        \item For the cases of three target (and five target) allocations, the recursive insertion construction method yields a lower minimum and mean deviation from the modified MD algorithm compared to the load balancing method by around $1$ to $5\%$. However, the maximum deviation of the solution from the recursive insertion method is marginally higher than that of the load-balancing method by around $1 \%$.
    \end{enumerate}
    It is now desired to pick the best construction method using the computation time of the heuristic as well. The distribution of the computation time for the considered variations in the neighborhood (except for the actual tour cost metric, which is not considered based on the discussion above) is shown in Fig.~\ref{fig: computation_time_box_plot_switch_swap}. It can be observed that the heuristic using the recursive insertion method typically has a lower maximum, mean, and median computation time than the heuristic using the load balancing method for a fixed metric choice and $n$ value
    by about $100$ to $200$ seconds, $20$ to $50$ seconds, and $10$ to $40$ seconds, respectively, 
    over all instances.
    
    From these observations, it can be concluded that the recursive insertion method yields a better solution at a lower computation time for the heuristic than the load balancing method and, therefore, would be the construction method considered henceforth.
    \item Using the selected recursive insertion method, for each value of $n$, it is desired to identify the metric to be chosen for the heuristic from the least insertion cost and least estimated tour cost metrics. Using the percentage deviation for each $n$ value given in  Fig.~\ref{subfig: variation_objective_val_zero_targ} for the zero-target case and Fig.~\ref{subfig: variation_objective_val_three_targ} for the three-target case, it can be observed that the least insertion cost metric produces a better solution. This is because the least insertion cost metric yields a lower maximum deviation, by at most around $2\%$, and minimum deviation, by around $1\%$, compared to the least estimated tour cost metric. 
    Additionally, from Fig.~\ref{fig: computation_time_box_plot_switch_swap}, it can be observed that both metrics yield similar computation time for the same $n$ value, with the least estimated tour cost metric yielding marginally better maximum and mean computation time by around $40$ seconds and $10$ seconds, respectively, for the zero-target cases for a fixed $n$ value. However, for the three-target and five-target cases, the least insertion cost metric performed faster than the least estimated tour cost metric in terms of the maximum and mean computation time by around $100$ seconds and $10$ seconds, respectively, for a fixed $n$ value. 
    
    From these observations, it can be observed that the least insertion cost metric is the best-performing metric for the switch and swap neighborhoods.
    \item Finally, it can be observed from Figs.~\ref{subfig: variation_objective_val_zero_targ} and \ref{subfig: variation_objective_val_three_targ} that with the initial solution obtained using the recursive insertion method with the least insertion cost metric, $n = 3$ yielded marginally better maximum deviation, mean deviation, and median deviation compared to $n = 2$ for the zero-target cases by around $0.1\%, 0.4\%,$ and $0.4\%,$ respectively. However, for the three-target (and five-target cases), $n = 2$ yielded better maximum and mean deviation compared to $n = 3$ by around $1$ to $2\%$ and $0.1\%,$ respectively. Noting additionally that the mean and maximum computation time increases from $n = 2$ to $n = 3$ by around $20$ to $40$ seconds and $100$ to $400$ seconds, respectively, as can be observed from Fig.~\ref{fig: computation_time_box_plot_switch_swap}, $n = 2$ is picked for the switch and swap neighborhoods.    
\end{enumerate}

Therefore, the recursive insertion method for constructing the initial feasible solution, and the least insertion cost metric with $n = 2$ for the switch and swap neighborhoods will be considered henceforth for the heuristic.

\textbf{Remark:} In Fig.~\ref{fig: computation_time_box_plot_switch_swap}, the distribution for the modified MD algorithm is not shown since the computation time for this algorithm is substantially higher than all the considered variations. 
This can be observed from the detailed computational results for the target switch and swap neighborhood analysis provided in Appendix~\ref{appsubsect: comp_results_switch_swap}.
Hence, to depict the difference in the distribution with different variations in the switch and swap neighborhoods, the computation time for the modified MD algorithm is not shown.


\subsection{Identifying sensitive instances} \label{subsect: sensitive_instances}

\begin{figure}[htb!]
    \centering
    \subfigure[Zero-target assignment instances \label{subfig: zero_target_assignment}]{\includegraphics[width = 0.48\textwidth]{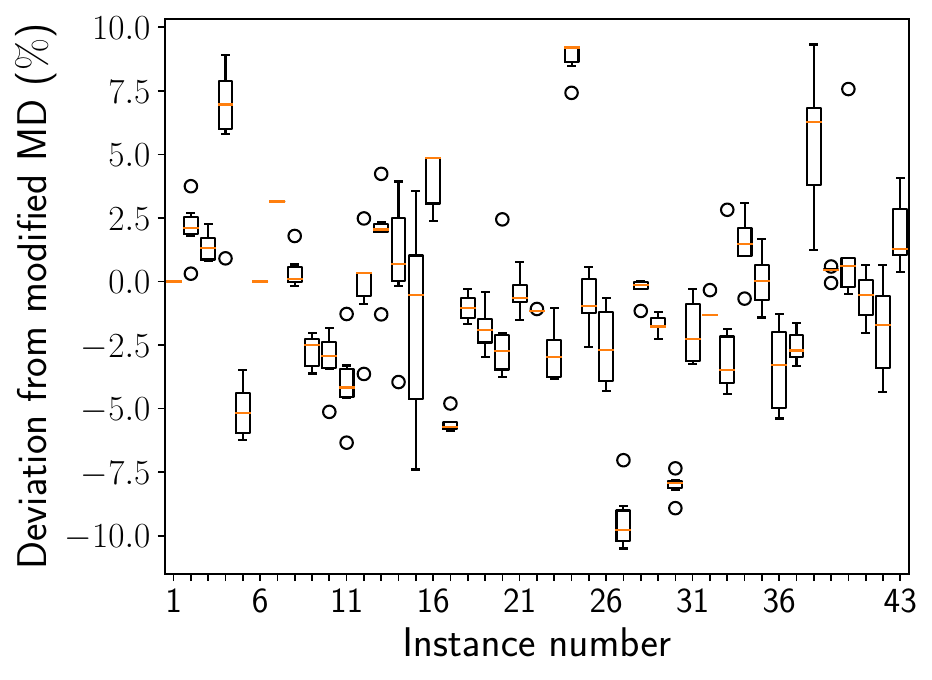}} \hfill
    \subfigure[Three-target assignment instances \label{subfig: three_target_assignment}]{\includegraphics[width = 0.48\textwidth]{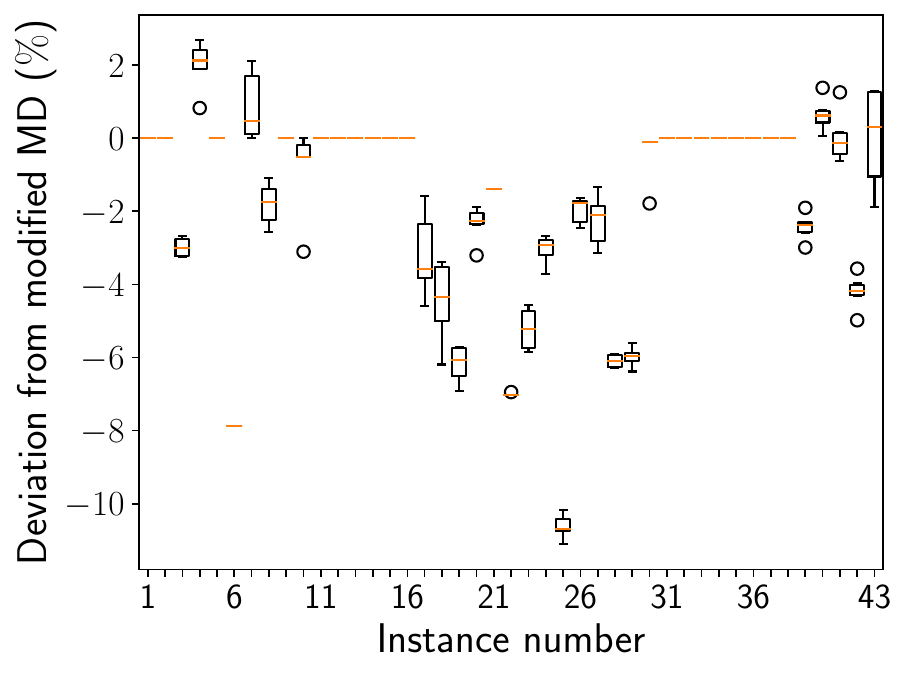}}
    \subfigure[Five-target assignment instances \label{subfig: five_target_assignment}]{\includegraphics[width = 0.48\textwidth]{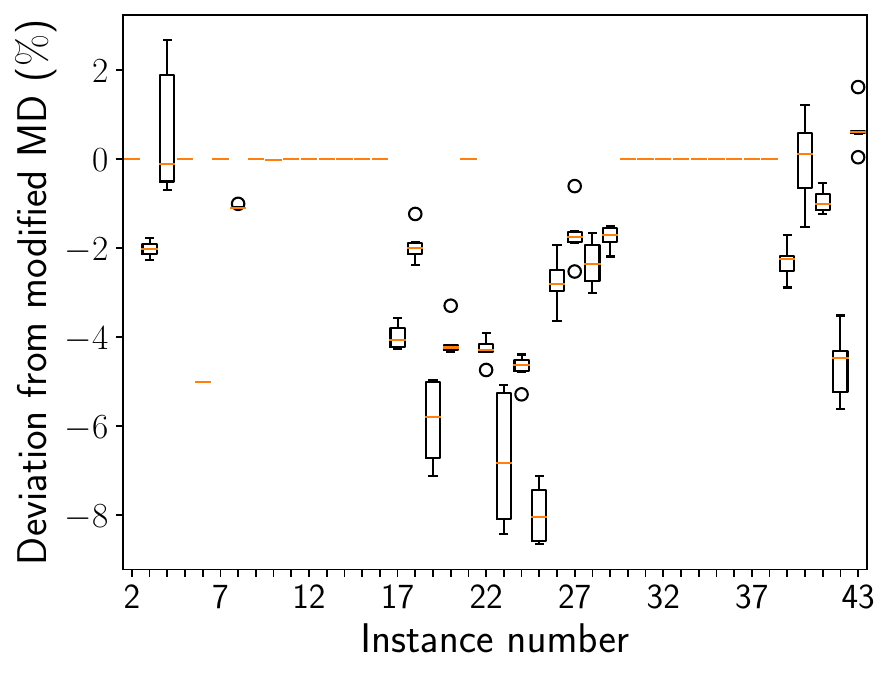}}
    \caption{Sensitivity of instances to change in metric and number of vehicles} \label{fig: sensitivity_scenarios_instances}
\end{figure}

From the previous analysis of the variations in the target switch and swap neighborhoods, the best-performing variation was selected for the considered neighborhoods.
It is now desired to utilize a multi-target swap neighborhood to improve the feasible solution obtained further. For this purpose, it is desired to identify a set of instances wherein the solution obtained is sensitive to variations in the parameters or metrics considered for the neighborhoods to simplify the analysis of the multi-target neighborhood to be used. To this end, the variation in the best feasible solution obtained for all instances for the heuristic with the target switch and swap neighborhoods utilizing
\begin{itemize}
    \item The recursive insertion method for constructing the initial feasible solution,
    \item The least insertion cost or the least estimated tour cost metrics for picking the vehicle for switching or swapping a target with and
    \item With $n = 1, 2, 3$,
\end{itemize}
are considered. The distribution in the solution obtained for the zero-target assignment, three-target assignment, and five-target assignment instances are shown in Figs.~\ref{subfig: zero_target_assignment}, \ref{subfig: three_target_assignment}, and \ref{subfig: five_target_assignment}, respectively. From these figures, it can be observed that
\begin{enumerate}
    \item The considered heuristic performs marginally worse than the modified MD algorithm on instances $4, 7, 24, 38, 42, 43$. 
    \item In addition, instances $18, 19, 23, 26, 27, 40$ can be observed to have large variations in the obtained feasible solution cost.
\end{enumerate}
Therefore, the analysis of multi-target swap neighborhoods will be restricted to the considered set of instances.

\subsection{Multi-Target Swap Neighborhood Analysis} \label{subsect: mult_targ_swap}

It must be recalled that two variations in the multi-target swap neighborhood were considered in Section~\ref{subsubsect: multi-target}: fixed-structure and variable-structure. In the fixed-structure multi-target swap neighborhood, two parameters were considered: $m,$ which denotes the size of target sets considered for swapping, and $n,$ which denotes the number of candidate solutions explored in the neighborhood corresponding to a selected target set. In this study, two variations in $m$ are considered, which are $m = 2,$ and $m = 2, 3,$ and two variations in $n$ are considered, which are $n = 10, 20.$ It should be noted that $m = 2, 3$ corresponds to utilizing two multi-target swap neighborhoods one after the other, wherein the first neighborhood corresponds to $m = 2$ and the second corresponds to $m = 3.$ For the variable-structure multi-target swap neighborhood, two parameters were similarly considered: $m$ and $n$. Since $m$ denotes the maximum target set size, including target sets of a smaller size as well, $m = 3$ and $4$ are considered variations, and $n = 10, 20$ are considered.

\begin{figure}[htb!]
    \centering
    \includegraphics[width=\linewidth]{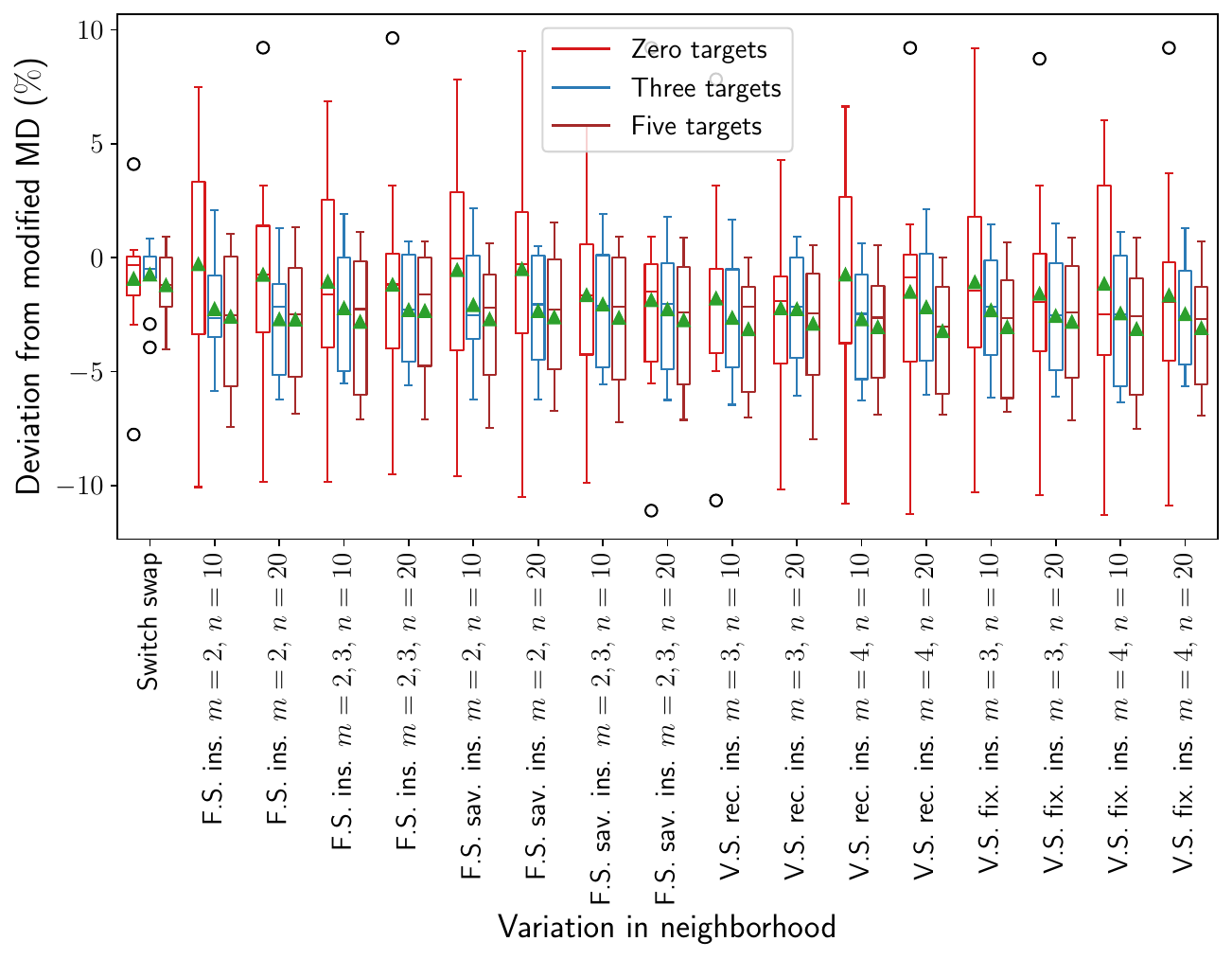}
    \caption{Box plot depicting the deviation of best solution with respect to the modified MD algorithm for each variation in the multi-target swap neighborhood for the three types of instances differing based on the number of targets allocated to each vehicle (along with the mean deviation depicted using a triangle)}
    \label{fig: deviation_box_plot}
\end{figure}

\begin{figure}[htb!]
    \centering
    \includegraphics[width=\linewidth]{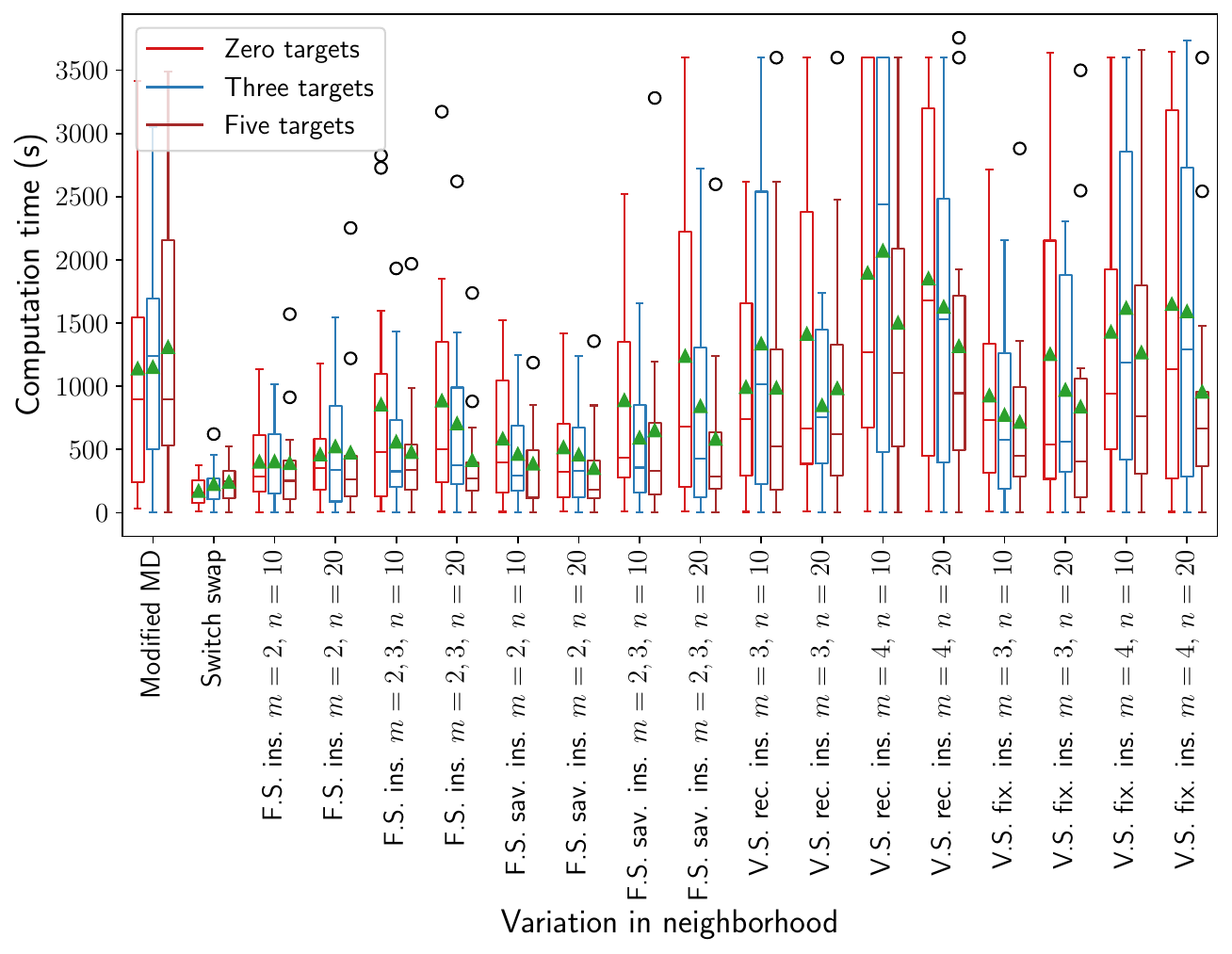}
    \caption{Box plot depicting the computation time for each variation in the multi-target swap neighborhood for the three types of instances differing based on the number of targets allocated to each vehicle (along with the mean time depicted using a triangle)}
    \label{fig: computation_time_box_plot}
\end{figure}

The depiction of the distribution of the deviation with respect to the modified MD algorithm and the computation time for all variations are shown in Figs.~\ref{fig: deviation_box_plot} and \ref{fig: computation_time_box_plot}, respectively, for the selected instances from the sensitivity study. Additionally, the deviation of the heuristic with the best target switch and swap neighborhoods from the modified MD algorithm and the distribution of the computation time are also shown in Figs.~\ref{fig: deviation_box_plot} and \ref{fig: computation_time_box_plot}, respectively.

\textbf{Remark:} The detailed computational results for this analysis are presented in Appendix~\ref{appsubsect: comp_results_mult_target}.


\begin{figure}[htb!]
    \centering
    \subfigure[Impact of variation in fixed-structure multi-target swap neighborhoods on deviation from modified MD algorithm for zero-target instances \label{subfig: variation_metric_obj_val_fixed_struct_zero_targ}]{\includegraphics[width = 0.48\textwidth]{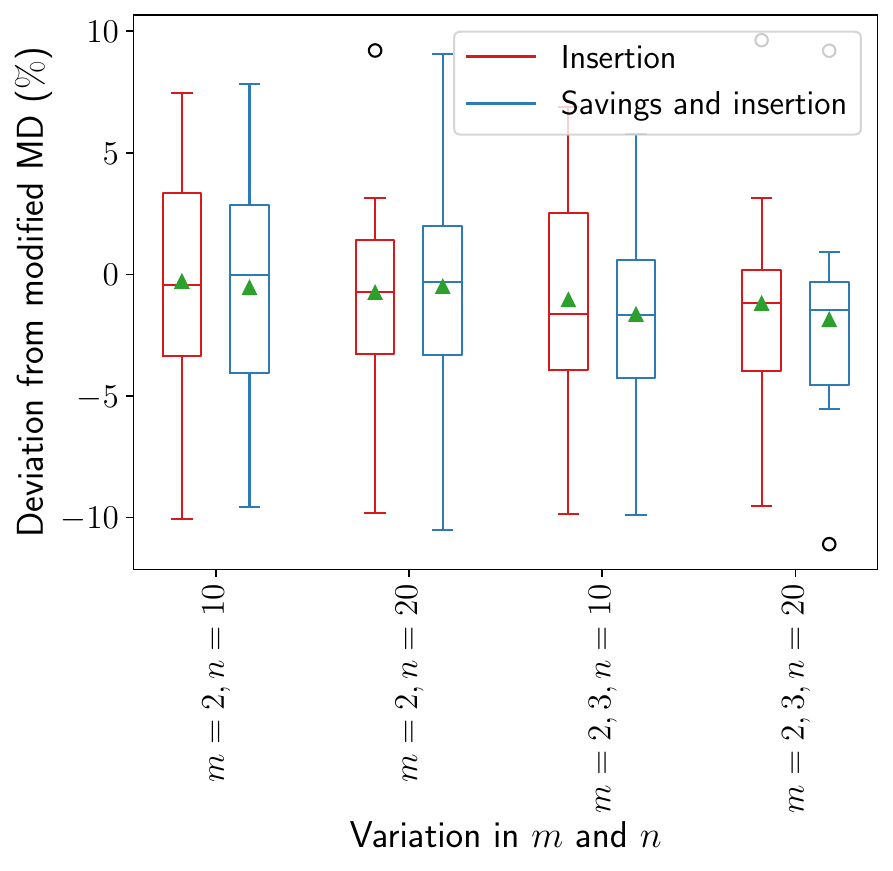}} \hfill
    \subfigure[Impact of variation in fixed-structure multi-target swap neighborhoods on deviation from modified MD algorithm for three-target instances \label{subfig: variation_metric_obj_val_fixed_struct_three_targ}]{\includegraphics[width = 0.48\textwidth]{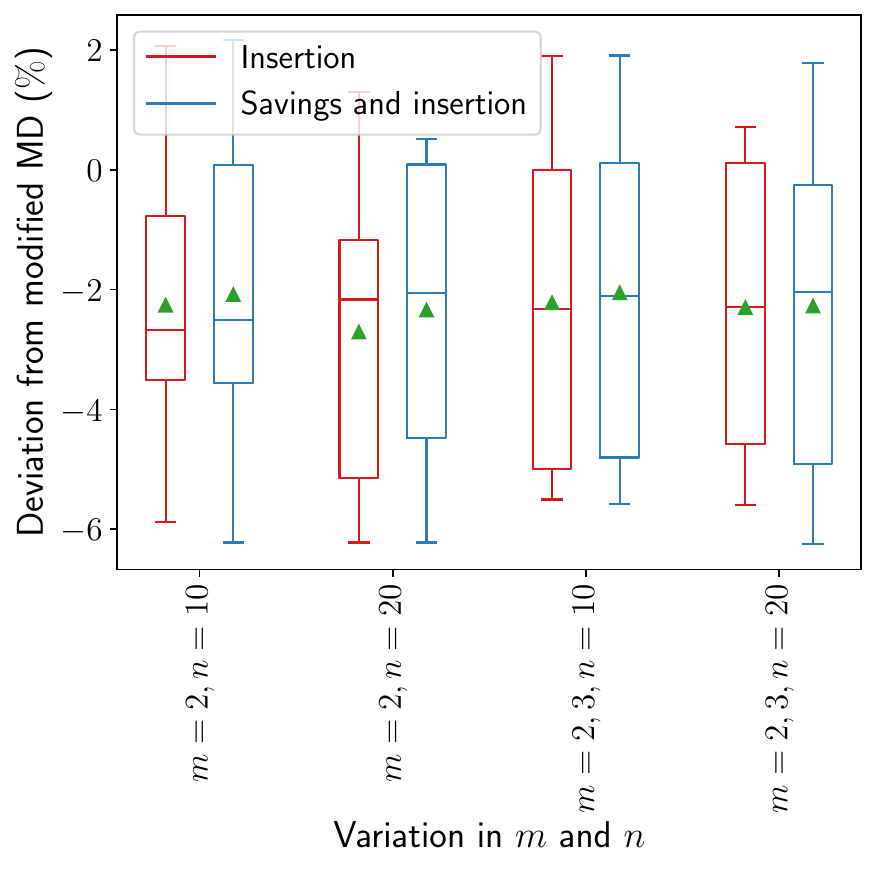}} \hfill
    \subfigure[Impact of variation in fixed-structure multi-target swap neighborhoods on computation time for zero-target instances \label{subfig: variation_metric_comp_time_fixed_struct_zero_targ}]{\includegraphics[width = 0.48\textwidth]{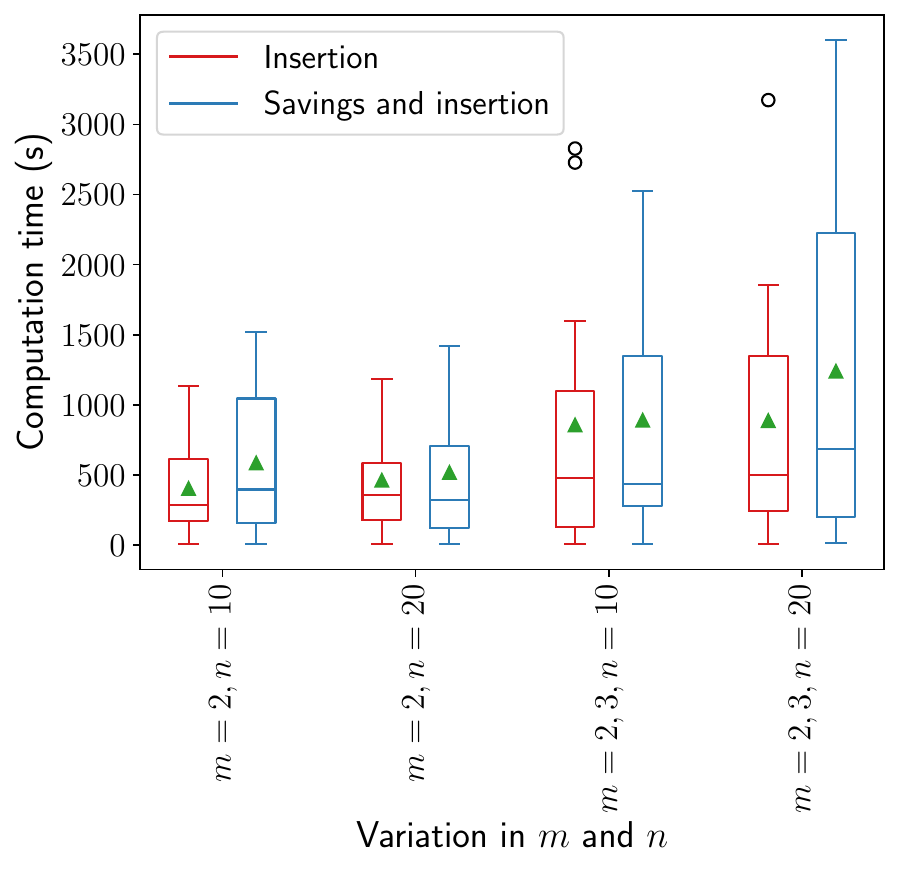}} \hfill
    \subfigure[Impact of variation in fixed-structure multi-target swap neighborhoods on computation time for three-target instances \label{subfig: variation_metric_comp_time_fixed_struct_three_targ}]{\includegraphics[width = 0.48\textwidth]{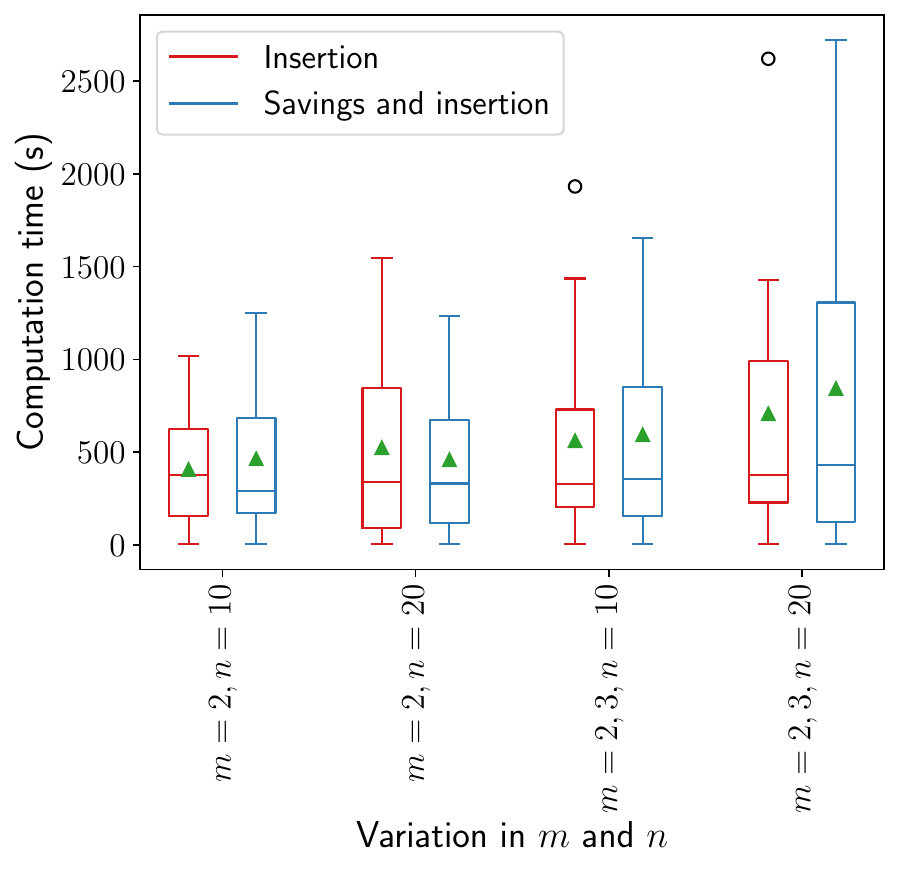}}
    \caption{Analysis of variations in fixed-structure multi-target neighborhood} \label{fig: analysis_fixed_struct_mult_targ_neighborhood}
\end{figure}

From Figs.~\ref{fig: deviation_box_plot} and \ref{fig: computation_time_box_plot}, the following observations can be made:
\begin{enumerate}
    \item Comparing the fixed-structure neighborhood variations and the variable-structure neighborhood variations, it can be observed that both variation types predominantly yield similar deviations in the objective value. The variable-structure neighborhood variations yield marginally better minimum, mean, and median deviation by about $1\%$ to $2\%$ for all three types of instances (based on the number of targets assigned to each vehicle). 
    However, the variable-structure neighborhood variations have a higher computation time. In particular, the mean and median computation times are higher for the variable-structure neighborhood variations than the corresponding fixed structure variations by $1.5$ to $3$ times, which is about $400$ to $1000$ seconds (depending on the parameters chosen). Therefore, the fixed-structure variations will be considered for the heuristic due to a substantially lower computation time at the cost of a marginal loss in the improvement of the objective value.
    \item The deviations for the variations in the fixed-structure multi-target swap neighborhood are shown in Figs.~\ref{subfig: variation_metric_obj_val_fixed_struct_zero_targ} and \ref{subfig: variation_metric_obj_val_fixed_struct_three_targ} for the zero-target and three-target cases, respectively, and the computation time for the same cases are shown in Figs.~\ref{subfig: variation_metric_comp_time_fixed_struct_zero_targ} and \ref{subfig: variation_metric_comp_time_fixed_struct_three_targ}, respectively. It should be noted that the figures for the five-target case are not shown due to their similarity to the three-target case. From Figs.~\ref{subfig: variation_metric_obj_val_fixed_struct_zero_targ} and \ref{subfig: variation_metric_obj_val_fixed_struct_three_targ}, it can be observed that the minimum, mean, and median deviations in the objective value are similar for all variations (except for $m = 2, 3$ with $n = 20$ for zero-target allocation instances, which showed a poor improvement from the modified MD algorithm). Therefore, the ideal metric choice will be picked based on the computation time.
    
    From Figs.~\ref{subfig: variation_metric_comp_time_fixed_struct_zero_targ} and \ref{subfig: variation_metric_comp_time_fixed_struct_three_targ}, it can be observed that the variations with the insertion cost metric perform better than the savings and insertion cost metric based on computation time, typically yielding a lower mean, median, and maximum computation time (except for $m = 2, n = 20$ for three-target and five-target cases) by around $50$ to $150$ seconds. Hence, the insertion cost metric for the fixed-structure multi-target swap neighborhood is the preferred metric to be utilized for the neighborhood.
    \item From Fig.~\ref{fig: computation_time_box_plot}, it can be observed that increasing $m$ from $2$ to $2, 3$ keeping $n$ constant for the chosen fixed-structure multi-target swap neighborhood with the insertion cost metric typically leads to an increase in the mean, median, and maximum computation time by about $1.5$ to $2$ times (except for $m = 2, n = 20$ for three-target and five-target cases). However, the increase in the minimum, mean, and median deviations when $m$ is increased is marginal and is at most about $1\%$ (refer to Fig.~\ref{fig: deviation_box_plot}). Hence, it is desired to pick $m = 2$ due to a lower computation time. 
    \item Since increasing $n$ from $10$ to $20$ for $m = 2$ leads to an increase in the mean and median computation time by around at most $50$ to $100$ seconds, a larger value of $n$ can be considered if an improvement in the objective value is observed. From Figs.~\ref{subfig: variation_metric_obj_val_fixed_struct_zero_targ} and \ref{subfig: variation_metric_obj_val_fixed_struct_three_targ}, it can be observed that increasing $n$ from $10$ to $20$ leads to a marginally improved mean and median deviations by about at most $0.5\%$ with marginal improvement in the minimum deviation for the three-target and five-target cases by about $0.5\%.$ Due to a marginal improvement in the deviations at the expense of a marginally increased computation time, $n = 20$ can be picked for the heuristic. 
\end{enumerate}

\subsection{Summary of results for heuristic}

Having selected the best heuristic in the previous step, a summary of the best solution obtained from the modified MD heuristic, the heuristic with the best switch and swap neighborhood, and the heuristic with the best switch and swap neighborhood along with the best multi-target swap neighborhood is provided for the $43$ heterogeneous instances with zero-target allocation, $43$ heterogeneous instances with three-target allocation, and $42$ heterogeneous instances with five-target instances in Tables~\ref{tab: summary_results_zero_target}, \ref{tab: summary_results_three_target}, and \ref{tab: summary_results_five_target}, respectively. The algorithm that provides the best solution for each instance is also highlighted in this table. From these tables, it can be observed that the heuristic with the switch and swap neighborhood and the multi-target swap neighborhood yielded the best solution for the highest number of instances, followed by the heuristic with the switch and swap neighborhood, and the modified MD algorithm. In particular,
\begin{enumerate}
    \item For the zero-target, three-target, and five-target instances, the heuristic utilizing the best multi-target swap neighborhood uniquely yielded the best solution on $17$, $12$, and $15$ instances, respectively, and provided a better solution than the modified MD algorithm on $28$, $22$, and $21$ instances, respectively. Furthermore, the heuristic with the multi-target swap yielded a solution with the same objective value as that of the modified MD algorithm on $2,$ $19,$ and $20$ instances for the zero-target, three-target, and five-target cases, respectively.
    \item For the zero-target, three-target, and five-target instances, the heuristic with the best switch and swap neighborhood uniquely yielded the best solution on $10,$ $8,$ and $4$ instances, respectively, and provided a better solution than the modified MD algorithm on $23,$ $20,$ and $21$ instances, respectively. Furthermore, the heuristic with the best switch and swap neighborhood yielded a solution with the same objective value as that of the modified MD algorithm on $2,$ $18,$ and $20$ instances for the zero-target, three-target, and five-target cases, respectively.
    \item In addition, both variations of the considered heuristic yield a solution faster than the modified MD algorithm, with the maximum computation time for the heuristic with the switch and swap neighborhood being about $820$ seconds and the heuristic with the switch and swap neighborhood and the multi-target swap neighborhood being about $2250$ seconds. It should be noted that for the heuristic with multi-target swap neighborhood, the computation time exceeded $1500$ seconds on only four out of the $128$ instances. On the other hand, for the heuristic with the switch and swap neighborhood, the computation time exceeded $600$ seconds on only one instance out of $128$ instances. In contrast, the maximum computation time for the modified MD algorithm was $3600$ seconds, with the algorithm exceeding $1500$ seconds on $27$ out of the $128$ instances.
\end{enumerate}
Hence, the heuristic with the best switch and swap neighborhood can be used when the user requires a quick, feasible solution, whereas the heuristic with the best switch and swap neighborhood and the best multi-target swap neighborhood can be used for a better feasible solution. 

\begin{figure}[htb!]
    \centering
    \subfigure[Modified MD algorithm \label{subfig: modified_MD_tour}]{\includegraphics[width = 0.6\textwidth]{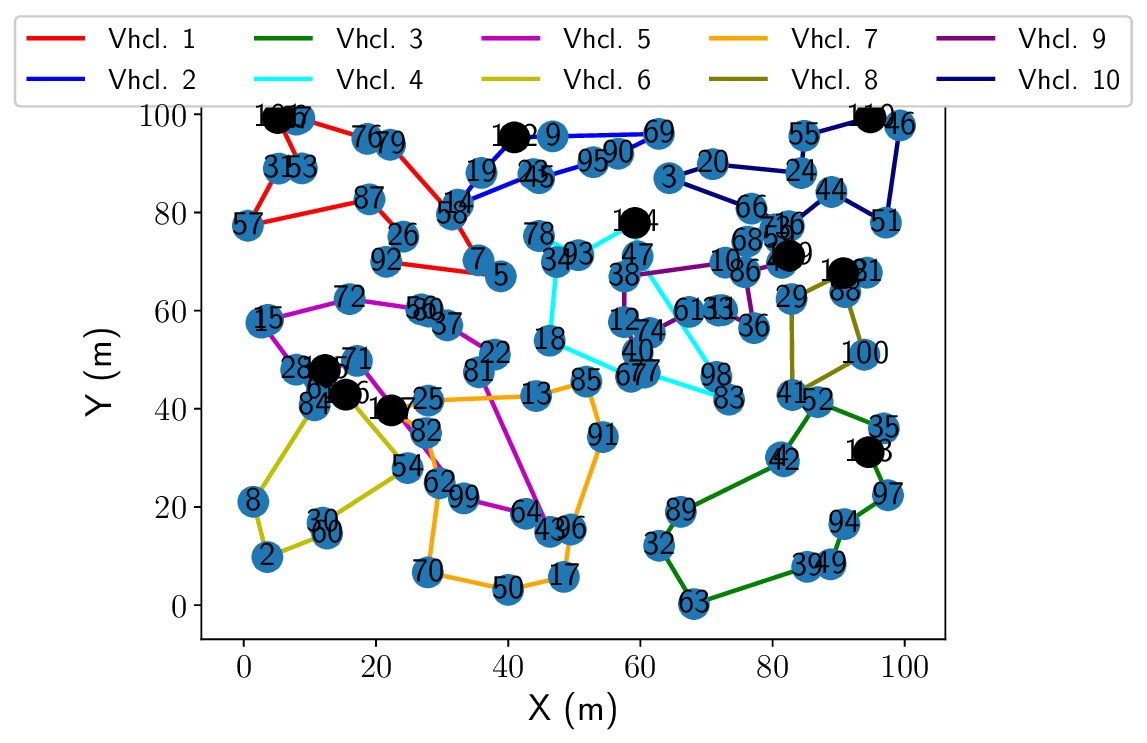}} \hfill
    \subfigure[Heuristic with best switch and swap neighborhood \label{subfig: heuristic_best_switch_swap_neighborhood}]{\includegraphics[width = 0.6\textwidth]{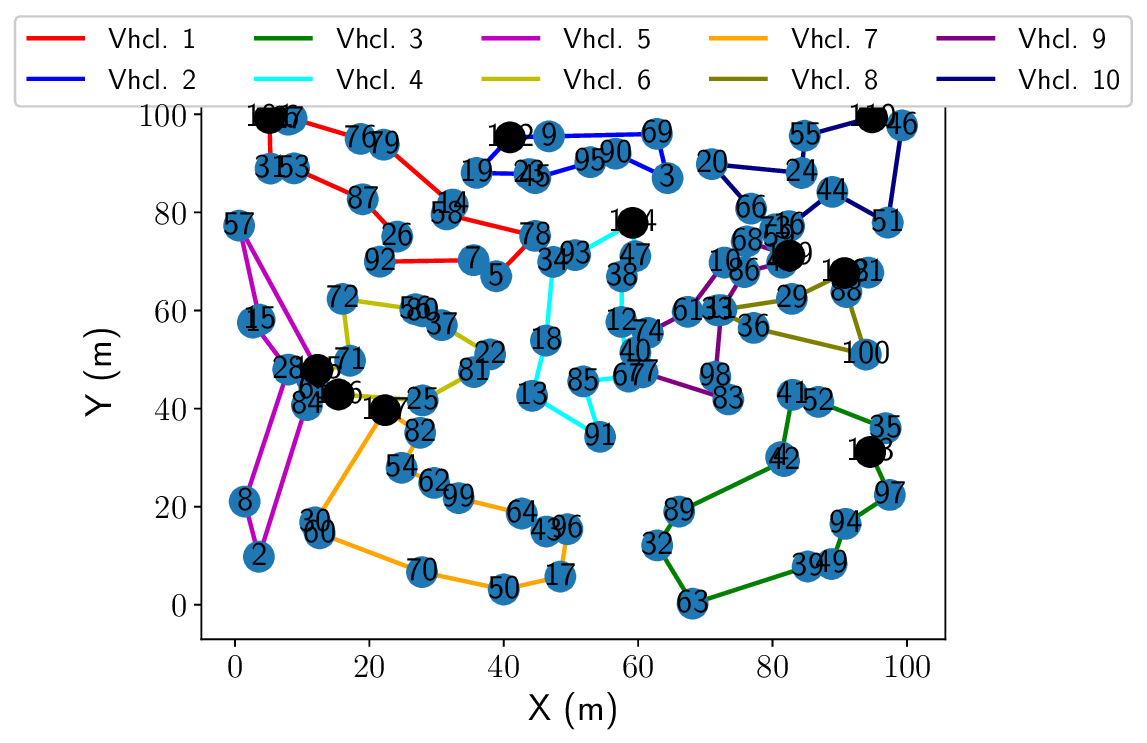}} \hfill
    \subfigure[Heuristic with best switch and swap neighborhood and best multi-target swap neighborhood \label{subfig: best_mult_targ_neighborhood_tours}]{\includegraphics[width = 0.6\textwidth]{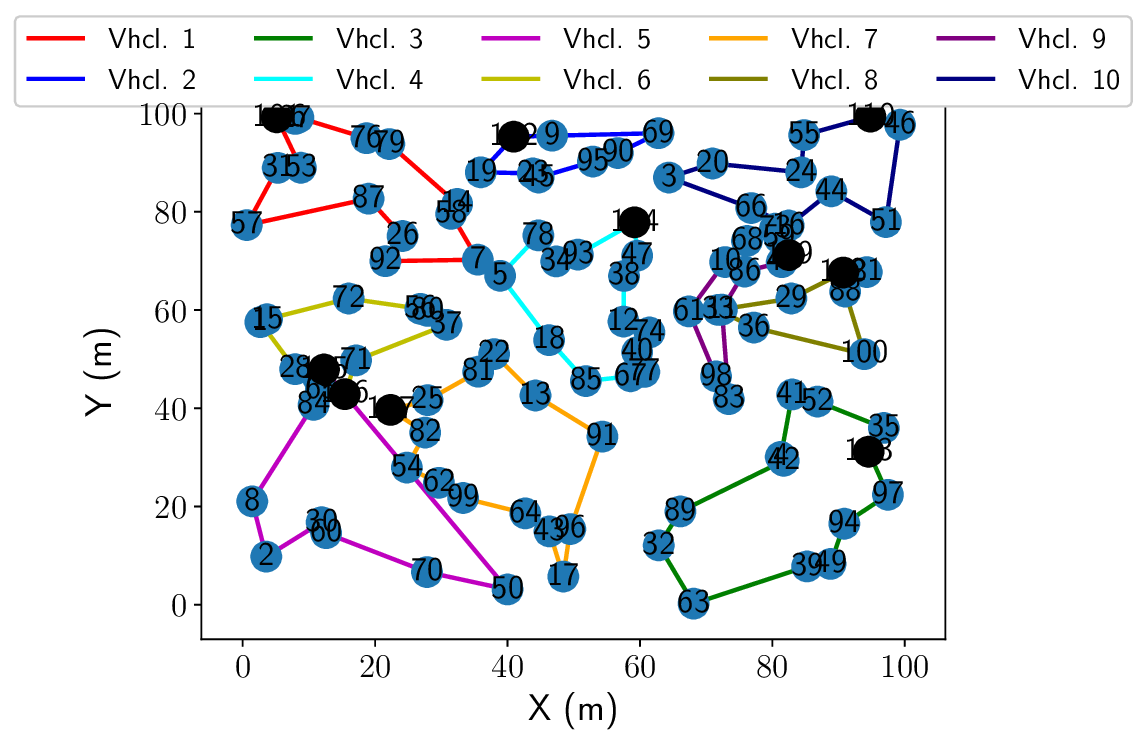}}
    \caption{Feasible tours generated for MM11 with no vehicle-target assignments} \label{fig: feasible_tours_generated}
\end{figure}

A demonstration of the tours obtained is shown in Fig.~\ref{fig: feasible_tours_generated} for instance MM11, which comprises $100$ targets and $10$ vehicles, for no vehicle-target allocations. In this instance, the objective value of the solution obtained from the modified MD algorithm, the heuristic with the best switch and swap neighborhood, and the heuristic with the best switch and swap neighborhood and best multi-target swap neighborhood are $73.49,$ $71.06,$ and $70.12,$ respectively.

\begin{table}[htb!]
    \centering
    \caption{Summary of results for heuristic with best switch and swap neighborhood, and best switch and swap neighborhood along with best multi-target swap neighborhood for zero target instances}
    \label{tab: summary_results_zero_target}
    \begin{tabular}{ccccccc}
    \toprule
    \multirow{2}{*}{\textbf{Instance}} & \multicolumn{3}{c}{\textbf{Best solution}} & \multicolumn{3}{c}{\textbf{Computation time (s)}} \\ \cmidrule{2-7} 
     & \begin{tabular}[c]{@{}c@{}}Modified \\ MD\\ algorithm\end{tabular} & \begin{tabular}[c]{@{}c@{}}Best\\ switch\\ \& swap\end{tabular} & \begin{tabular}[c]{@{}c@{}}Best switch \&\\ swap, and best\\ mult. targ. swap\end{tabular} & \begin{tabular}[c]{@{}c@{}}Modified \\ MD\\ algorithm\end{tabular} & \begin{tabular}[c]{@{}c@{}}Best\\ switch\\ \& swap\end{tabular} & \begin{tabular}[c]{@{}c@{}}Best switch \& \\ swap, and best\\ mult. targ. swap\end{tabular} \\ \midrule
    MM1\_0 & \textbf{136.72} & \textbf{136.72} & \textbf{136.72} & 6.60 & 3.76 & 5.03 \\
    MM2\_0 & \textbf{84.88} & 86.39 & 86.39 & 68.50 & 40.17 & 77.03 \\
    MM3\_0 & 170.21 & 171.61 & \textbf{169.49} & 98.32 & 67.74 & 70.56 \\
    MM4\_0 & 302.52 & 305.27 & \textbf{302.21} & 1000.44 & 372.78 & 1181.45 \\
    MM5\_0 & 215.48 & \textbf{203.58} & 206.99 & 273.52 & 180.63 & 168.28 \\
    MM6\_0 & \textbf{65.68} & \textbf{65.68} & \textbf{65.68} & 460.28 & 54.90 & 150.68 \\
    MM7\_0 & \textbf{126.81} & 130.80 & 130.80 & 30.59 & 4.22 & 5.71 \\
    MM8\_0 & 131.90 & 132.78 & \textbf{131.14} & 353.65 & 27.65 & 62.13 \\
    MM9\_0 & 112.68 & 110.15 & \textbf{109.51} & 185.69 & 110.85 & 299.55 \\
    MM10\_0 & 118.91 & 116.12 & \textbf{111.18} & 24.30 & 13.07 & 19.99 \\
    MM11\_0 & 73.49 & 71.06 & \textbf{70.12} & 45.59 & 21.46 & 18.22 \\
    MM12\_0 & 50.02 & 51.26 & \textbf{49.51} & 17.51 & 20.40 & 17.96 \\
    MM13\_0 & 82.67 & 84.61 & \textbf{81.46} & 79.61 & 38.59 & 79.04 \\
    MM14\_0 & 77.84 & 80.90 & \textbf{75.41} & 100.05 & 27.08 & 79.98 \\
    MM15\_0 & \textbf{76.52} & 76.86 & 77.45 & 53.09 & 51.12 & 22.25 \\
    MM16\_0 & \textbf{61.61} & 63.13 & 67.12 & 269.64 & 67.13 & 80.31 \\
    MM17\_0 & 173.09 & 164.78 & \textbf{162.66} & 1072.74 & 257.04 & 624.49 \\
    MM18\_0 & 239.15 & 237.81 & \textbf{235.83} & 2937.39 & 258.22 & 1183.31 \\
    MM19\_0 & 253.66 & 247.48 & \textbf{247.35} & 1520.33 & 243.05 & 427.39 \\
    MM20\_0 & 201.44 & \textbf{193.89} & 198.72 & 2002.70 & 229.15 & 308.56 \\
    MM21\_0 & 147.05 & \textbf{144.83} & 145.39 & 192.37 & 41.73 & 100.96 \\
    MM22\_0 & 241.09 & \textbf{238.27} & \textbf{238.27} & 783.91 & 43.35 & 116.11 \\
    MM23\_0 & 330.00 & 318.24 & \textbf{317.32} & 111.48 & 51.75 & 106.84 \\
    MM24\_0 & \textbf{111.06} & 121.28 & 121.29 & 175.14 & 32.72 & 71.21 \\
    MM25\_0 & 189.84 & \textbf{184.94} & 185.56 & 536.81 & 45.26 & 166.47 \\
    MM26\_0 & 270.59 & 259.48 & \textbf{257.79} & 1431.36 & 111.38 & 434.45 \\
    MM27\_0 & 245.31 & \textbf{220.69} & 221.19 & 384.62 & 70.79 & 210.53 \\
    MM28\_0 & 238.73 & 238.65 & \textbf{235.22} & 2132.53 & 493.84 & 611.16 \\
    MM29\_0 & 301.84 & \textbf{296.37} & 297.41 & 1890.68 & 361.26 & 455.36 \\
    MM30\_0 & 102.16 & \textbf{94.05} & 94.60 & 197.41 & 57.07 & 65.91 \\
    MM31\_0 & 82.91 & 82.66 & \textbf{80.59} & 91.60 & 40.55 & 66.88 \\
    MM32\_0 & 55.69 & \textbf{54.95} & \textbf{54.95} & 19.07 & 4.58 & 7.06 \\
    MM33\_0 & 56.74 & \textbf{54.22} & 55.39 & 254.22 & 54.98 & 53.8 \\
    MM34\_0 & \textbf{59.33} & 59.94 & 62.06 & 161.51 & 93.34 & 150.88 \\
    MM35\_0 & \textbf{64.26} & 64.72 & 65.07 & 104.07 & 49.47 & 74.76 \\
    MM36\_0 & 94.99 & 93.79 & \textbf{92.86} & 809.00 & 197.61 & 472.27 \\
    MM37\_0 & 80.97 & \textbf{79.42} & 84.52 & 111.64 & 82.83 & 39.64 \\
    MM38\_0 & \textbf{106.75} & 116.70 & 107.91 & 261.26 & 153.19 & 205.39 \\
    MM39\_0 & \textbf{132.48} & 133.07 & 133.07 & 403.05 & 31.04 & 62.26 \\
    MM40\_0 & \textbf{193.69} & 195.45 & 197.06 & 786.78 & 208.39 & 281.48 \\
    MM41\_0 & \textbf{174.28} & 175.43 & 174.46 & 3400.71 & 400.78 & 543.88 \\
    MM42\_0 & 243.94 & \textbf{234.77} & 236.35 & 1632.53 & 823.75 & 964.83 \\
    MM43\_0 & \textbf{226.11} & 235.28 & 229.01 & 3415.71 & 157.56 & 457.07 \\ \bottomrule
    \end{tabular}
\end{table}

\begin{table}[htb!]
    \centering
    \caption{Summary of results for heuristic with best switch and swap neighborhood, and best switch and swap neighborhood along with best multi-target swap neighborhood for three target instances}
    \label{tab: summary_results_three_target}
    \begin{tabular}{ccccccc}
    \toprule
    \multirow{2}{*}{\textbf{Instance}} & \multicolumn{3}{c}{\textbf{Best solution}} & \multicolumn{3}{c}{\textbf{Computation time (s)}} \\ \cmidrule{2-7} 
     & \begin{tabular}[c]{@{}c@{}}Modified \\ MD\\ algorithm\end{tabular} & \begin{tabular}[c]{@{}c@{}}Best\\ switch\\ \& swap\end{tabular} & \begin{tabular}[c]{@{}c@{}}Best switch \&\\ swap, and best\\ mult. targ. swap\end{tabular} & \begin{tabular}[c]{@{}c@{}}Modified \\ MD\\ algorithm\end{tabular} & \begin{tabular}[c]{@{}c@{}}Best\\ switch\\ \& swap\end{tabular} & \begin{tabular}[c]{@{}c@{}}Best switch \& \\ swap, and best\\ mult. targ. swap\end{tabular} \\ \midrule
    MM1\_3 & \textbf{209.41} & \textbf{209.41} & \textbf{209.41} & 1.13 & 1.20 & 1.28 \\
    MM2\_3 & \textbf{204.80} & \textbf{204.80} & \textbf{204.80} & 8.71 & 5.29 & 5.50 \\
    MM3\_3 & 232.38 & 225.02 & \textbf{224.73} & 268.47 & 53.33 & 51.94 \\
    MM4\_3 & \textbf{391.09} & 394.32 & 396.19 & 1648.87 & 414.07 & 698.55 \\
    MM5\_3 & \textbf{405.76} & \textbf{405.76} & \textbf{405.76} & 25.19 & 12.09 & 11.34 \\
    MM6\_3 & 175.01 & \textbf{161.24} & \textbf{161.24} & 73.49 & 34.56 & 48.51 \\
    MM7\_3 & \textbf{171.48} & 172.28 & 172.28 & 5.08 & 6.28 & 4.06 \\
    MM8\_3 & 194.62 & \textbf{189.93} & 191.25 & 289.56 & 66.17 & 104.34 \\
    MM9\_3 & \textbf{182.94} & \textbf{182.94} & \textbf{182.94} & 63.60 & 11.87 & 12.45 \\
    MM10\_3 & 165.55 & \textbf{164.71} & 165.55 & 13.20 & 10.15 & 7.31 \\
    MM11\_3 & \textbf{237.56} & \textbf{237.56} & \textbf{237.56} & 3.95 & 4.68 & 4.53 \\
    MM12\_3 & \textbf{251.92} & \textbf{251.92} & \textbf{251.92} & 5.67 & 6.62 & 6.60 \\
    MM13\_3 & \textbf{226.03} & \textbf{226.03} & \textbf{226.03} & 4.65 & 4.90 & 4.83 \\
    MM14\_3 & \textbf{190.10} & \textbf{190.10} & \textbf{190.10} & 7.39 & 5.49 & 5.42 \\
    MM15\_3 & \textbf{255.39} & \textbf{255.39} & \textbf{255.39} & 7.64 & 7.51 & 7.46 \\
    MM16\_3 & \textbf{268.74} & \textbf{268.74} & \textbf{268.74} & 20.55 & 15.22 & 16.06 \\
    MM17\_3 & 179.29 & 175.67 & \textbf{170.69} & 2229.42 & 341.62 & 825.01 \\
    MM18\_3 & 275.68 & \textbf{261.69} & 262.33 & 1277.53 & 476.91 & 1220.95 \\
    MM19\_3 & 288.40 & 271.86 & \textbf{271.00} & 1197.00 & 357.12 & 763.45 \\
    MM20\_3 & 231.48 & \textbf{224.06} & 224.42 & 3281.76 & 479.02 & 1589.29 \\
    MM21\_3 & 235.56 & \textbf{232.27} & \textbf{232.27} & 10.86 & 4.99 & 5.31 \\
    MM22\_3 & 274.27 & \textbf{254.98} & 255.46 & 1677.48 & 58.90 & 270.38 \\
    MM23\_3 & 376.40 & 355.76 & \textbf{352.99} & 538.92 & 81.22 & 101.58 \\
    MM24\_3 & 168.46 & \textbf{162.22} & 162.92 & 594.23 & 283.00 & 132.15 \\
    MM25\_3 & 249.19 & \textbf{221.51} & 224.06 & 500.84 & 179.52 & 326.84 \\
    MM26\_3 & 301.83 & 296.87 & \textbf{295.23} & 1292.75 & 194.81 & 133.14 \\
    MM27\_3 & 275.54 & \textbf{269.63} & 269.64 & 397.59 & 103.20 & 50.10 \\
    MM28\_3 & 274.06 & 257.90 & \textbf{254.35} & 2067.86 & 400.68 & 1211.76 \\
    MM29\_3 & 343.90 & 323.77 & \textbf{322.86} & 1855.56 & 247.67 & 438.41 \\
    MM30\_3 & 142.29 & \textbf{142.16} & \textbf{142.16} & 142.47 & 22.79 & 22.95 \\
    MM31\_3 & \textbf{204.49} & \textbf{204.49} & \textbf{204.49} & 7.99 & 5.73 & 6.15 \\
    MM32\_3 & \textbf{124.75} & \textbf{124.75} & \textbf{124.75} & 2.31 & 2.25 & 2.44 \\
    MM33\_3 & \textbf{184.10} & \textbf{184.10} & \textbf{184.10} & 9.97 & 7.61 & 8.18 \\
    MM34\_3 & \textbf{222.55} & \textbf{222.55} & \textbf{222.55} & 24.03 & 18.93 & 19.40 \\
    MM35\_3 & \textbf{185.46} & \textbf{185.46} & \textbf{185.46} & 5.18 & 5.03 & 5.30 \\
    MM36\_3 & \textbf{176.38} & \textbf{176.38} & \textbf{176.38} & 17.22 & 12.87 & 12.60 \\
    MM37\_3 & \textbf{237.96} & \textbf{237.96} & \textbf{237.96} & 14.95 & 11.59 & 12.18 \\
    MM38\_3 & \textbf{166.18} & \textbf{166.18} & \textbf{166.18} & 37.07 & 7.21 & 7.17 \\
    MM39\_3 & 163.88 & 159.61 & \textbf{158.66} & 635.10 & 145.45 & 206.07 \\
    MM40\_3 & 237.75 & 241.03 & \textbf{234.04} & 1945.05 & 78.97 & 540.32 \\
    MM41\_3 & 203.12 & 203.44 & \textbf{200.44} & 1970.77 & 342.38 & 511.04 \\
    MM42\_3 & 277.46 & 265.52 & \textbf{260.37} & 3050.48 & 519.84 & 1546.86 \\
    MM43\_3 & 244.93 & 248.07 & \textbf{240.40} & 1837.49 & 214.85 & 1084.33 \\ \bottomrule
    \end{tabular}
\end{table}

\begin{table}[htb!]
    \centering
    \caption{Summary of results for heuristic with best switch and swap neighborhood, and best switch and swap neighborhood along with best multi-target swap neighborhood for five target instances}
    \label{tab: summary_results_five_target}
    \begin{tabular}{ccccccc}
    \toprule
    \multirow{2}{*}{\textbf{Instance}} & \multicolumn{3}{c}{\textbf{Best solution}} & \multicolumn{3}{c}{\textbf{Computation time (s)}} \\ \cmidrule{2-7} 
     & \begin{tabular}[c]{@{}c@{}}Modified \\ MD\\ algorithm\end{tabular} & \begin{tabular}[c]{@{}c@{}}Best\\ switch\\ \& swap\end{tabular} & \begin{tabular}[c]{@{}c@{}}Best switch \&\\ swap, and best\\ mult. targ. swap\end{tabular} & \begin{tabular}[c]{@{}c@{}}Modified \\ MD\\ algorithm\end{tabular} & \begin{tabular}[c]{@{}c@{}}Best\\ switch\\ \& swap\end{tabular} & \begin{tabular}[c]{@{}c@{}}Best switch \& \\ swap, and best\\ mult. targ. swap\end{tabular} \\ \midrule
    MM2\_5 & \textbf{209.35} & \textbf{209.35} & \textbf{209.35} & 12.14 & 5.28 & 5.56 \\
    MM3\_5 & 232.06 & 226.77 & \textbf{226.62} & 310.90 & 29.35 & 37.40 \\
    MM4\_5 & 389.59 & 387.42 & \textbf{386.74} & 1591.12 & 294.76 & 450.63 \\
    MM5\_5 & \textbf{402.35} & \textbf{402.35} & \textbf{402.35} & 33.82 & 12.78 & 12.60 \\
    MM6\_5 & 180.39 & \textbf{171.35} & \textbf{171.35} & 165.38 & 20.30 & 22.78 \\
    MM7\_5 & \textbf{212.54} & \textbf{212.54} & \textbf{212.54} & 1.29 & 1.31 & 1.22 \\
    MM8\_5 & 192.35 & \textbf{190.23} & 190.24 & 155.77 & 15.09 & 23.61 \\
    MM9\_5 & \textbf{188.48} & \textbf{188.48} & \textbf{188.48} & 52.47 & 12.11 & 13.96 \\
    MM10\_5 & 161.83 & \textbf{161.80} & \textbf{161.80} & 11.69 & 3.92 & 4.16 \\
    MM11\_5 & \textbf{267.33} & \textbf{267.33} & \textbf{267.33} & 4.10 & 4.61 & 4.64 \\
    MM12\_5 & \textbf{242.07} & \textbf{242.07} & \textbf{242.07} & 5.52 & 6.51 & 6.54 \\
    MM13\_5 & \textbf{259.23} & \textbf{259.23} & \textbf{259.23} & 5.43 & 4.83 & 4.77 \\
    MM14\_5 & \textbf{254.36} & \textbf{254.36} & \textbf{254.36} & 5.26 & 5.47 & 5.42 \\
    MM15\_5 & \textbf{277.17} & \textbf{277.17} & \textbf{277.17} & 9.63 & 7.23 & 7.61 \\
    MM16\_5 & \textbf{250.96} & \textbf{250.96} & \textbf{250.96} & 20.77 & 14.52 & 15.18 \\
    MM17\_5 & 194.53 & 187.59 & \textbf{186.04} & 1672.29 & 79.51 & 792.63 \\
    MM18\_5 & 275.10 & 268.52 & \textbf{265.96} & 2339.94 & 409.52 & 1220.74 \\
    MM19\_5 & 292.97 & 273.12 & \textbf{272.94} & 953.90 & 284.92 & 456.62 \\
    MM20\_5 & 238.21 & 228.02 & \textbf{224.71} & 1042.59 & 416.37 & 2065.26 \\
    MM21\_5 & \textbf{239.06} & \textbf{239.06} & \textbf{239.06} & 7.00 & 4.46 & 4.58 \\
    MM22\_5 & 273.51 & 262.24 & \textbf{261.76} & 1611.01 & 184.44 & 218.30 \\
    MM23\_5 & 372.71 & 353.60 & \textbf{351.50} & 837.83 & 28.23 & 137.87 \\
    MM24\_5 & 167.64 & \textbf{158.78} & 159.10 & 722.31 & 178.55 & 97.15 \\
    MM25\_5 & 243.61 & 223.02 & \textbf{222.06} & 547.26 & 135.91 & 231.42 \\
    MM26\_5 & 311.97 & 302.91 & \textbf{302.10} & 567.54 & 256.70 & 405.23 \\
    MM27\_5 & 274.92 & 270.46 & \textbf{269.92} & 412.35 & 95.61 & 141.77 \\
    MM28\_5 & 267.52 & 263.06 & \textbf{260.05} & 1971.69 & 439.40 & 792.82 \\
    MM29\_5 & 353.32 & \textbf{346.69} & 348.14 & 3600.00 & 126.81 & 168.65 \\
    MM30\_5 & \textbf{162.11} & \textbf{162.11} & \textbf{162.11} & 52.50 & 11.82 & 12.06 \\
    MM31\_5 & \textbf{230.06} & \textbf{230.06} & \textbf{230.06} & 5.76 & 5.62 & 6.60 \\
    MM32\_5 & \textbf{118.07} & \textbf{118.07} & \textbf{118.07} & 2.08 & 2.25 & 2.35 \\
    MM33\_5 & \textbf{183.88} & \textbf{183.88} & \textbf{183.88} & 10.86 & 7.53 & 7.80 \\
    MM34\_5 & \textbf{218.08} & \textbf{218.08} & \textbf{218.08} & 26.29 & 17.26 & 18.30 \\
    MM35\_5 & \textbf{195.28} & \textbf{195.28} & \textbf{195.28} & 4.71 & 4.99 & 5.34 \\
    MM36\_5 & \textbf{187.81} & \textbf{187.81} & \textbf{187.81} & 19.25 & 11.67 & 11.91 \\
    MM37\_5 & \textbf{259.23} & \textbf{259.23} & \textbf{259.23} & 12.08 & 11.06 & 11.34 \\
    MM38\_5 & \textbf{168.39} & \textbf{168.39} & \textbf{168.39} & 50.42 & 7.48 & 6.19 \\
    MM39\_5 & 159.64 & \textbf{155.03} & 157.70 & 444.02 & 153.97 & 99.09 \\
    MM40\_5 & 241.33 & 240.70 & \textbf{239.89} & 2095.01 & 249.33 & 305.47 \\
    MM41\_5 & 201.27 & 199.12 & \textbf{198.58} & 1097.63 & 190.30 & 252.77 \\
    MM42\_5 & 274.82 & 262.28 & \textbf{256.37} & 3489.23 & 339.24 & 2252.52 \\
    MM43\_5 & \textbf{241.93} & 243.29 & 245.13 & 2656.53 & 324.10 & 227.77 \\ \bottomrule
    \end{tabular}
\end{table}

\section{Conclusions \& Remarks}

In this paper, a computationally efficient three-staged heuristic for generating high-quality feasible solutions for a heterogeneous min-max multi-vehicle multi-depot Traveling Salesman Problem was systematically developed. In this study, vehicles were considered to be functionally heterogeneous due to different vehicle-target assignments originating from different sensing capabilities of vehicles and structurally heterogeneous due to different (but not necessarily distinct) speeds. A total of $128$ instances were developed with varying numbers of vehicle-target assignments, varying speeds, varying distributions from which the targets are generated, and varying target-to-vehicle ratios. For the proposed heuristic, two variations in the generation of the initial feasible solution in the first stage were considered, and three neighborhoods were considered for the local search in the second stage. Multiple variations in the neighborhoods were considered, including the metric used to pick targets to be removed or swapped, the number of candidate vehicles to be considered for an exchange, and the number of candidate solutions explored in the neighborhood. Each variation was benchmarked against an existing heuristic in the literature, and improvements in terms of the objective value and computation time were recorded. Two variations in the heuristic were finally considered, and one of the two heuristics can be utilized based on the importance provided to the objective value or computational efficiency. In particular, it was observed that the heuristic that was selected for the best objective value outperforms the benchmarking heuristic on $71$ out of $128$ instances and provides a similar objective value on $41$ instances with a lower computation time. On the other hand, the heuristic selected based on the computation time produced a feasible solution within 10 minutes on $127$ out of the $128$ instances. Further, this heuristic variation produced a solution with a better objective value on $64$ out of $128$ instances than the benchmarking heuristic and produced a solution with the same objective value on $40$ instances as that of the benchmarking heuristic.

\backmatter


\bmhead{Acknowledgements}

The authors thank Dr. Kaarthik Sundar for providing a baseline code for implementing a Variable Neighborhood Search and for useful discussions pertaining to the development of the heuristic. \\

This work was cleared for public release under AFRL-2024-6089 and the views or opinions expressed in this work are those of the authors and do not reflect any position or opinion of the United States Government, US Air Force, or Air Force Research Laboratory.

\section*{Declarations}

\begin{itemize}
\item Conflict of interest/Competing interests: Not applicable
\item Ethics approval and consent to participate: Not applicable
\item Consent for publication: Not applicable
\item Consent to participate: Not applicable
\item Ethics approval: Not applicable
\item Author contribution:
\begin{itemize}
    \item Deepak Prakash Kumar: Implementation, Heuristic design, data analysis, manuscript writing
    \item Sivakumar Rathinam: Study conception and design, editing
    \item Swaroop Darbha: Study conception and design, editing
    \item Trevor Bihl: Editing, problem formulation
\end{itemize}
\item Data availability: The data can be made available for replication at the request of a user.
\end{itemize}

\begin{appendices}

\section{Detailed computational results}

In this section, detailed computational results obtained for the analysis of variations in the target switch and swap neighborhoods, and the computational results for the variations in the multi-target swap neighborhood are provided.

\subsection{Switch and swap neighborhood analysis} \label{appsubsect: comp_results_switch_swap}

The detailed results for the target switch and swap neighborhood analysis described in Section~\ref{subsubsect: switch_swap_analysis} are provided in this section. In particular, the description of the minimum, maximum, mean, and median deviation of the objective value with respect to the modified MD algorithm and the minimum, maximum, mean, and median computation time are provided in Tables~\ref{tab: no_target_allocation_soln_comparison_switch_swap}, \ref{tab: three_target_allocation_soln_comparison_switch_swap}, and \ref{tab: five_target_allocation_soln_comparison_switch_swap}, for the zero-target, three-target, and five-target instances, respectively.

\begin{sidewaystable}[htb!]
    \centering
    \caption{Comparison of deviation and computation time for zero-target allocation instances for the variations in the switch and swap neighborhoods}
    \label{tab: no_target_allocation_soln_comparison_switch_swap}
    \begin{tabular}{ccccccccccc}
    \toprule
    \multirow{2}{*}{\textbf{Metric}} & \multirow{2}{*}{$n$} & \multirow{2}{*}{\textbf{Construction}} & \multicolumn{4}{c}{\textbf{Deviation from modified MD ($\%$)}} & \multicolumn{4}{c}{\textbf{Computation time (s)}} \\ \cmidrule{4-11} 
     &  &  & Min & Max & Mean & Median & Min & Max & Mean & Median \\ \midrule
    \begin{tabular}[c]{@{}c@{}}Modified\\ MD\end{tabular} & - & Load balancing & - & - & - & - & 6.60 & 3415.71 & 695.07 & 261.26 \\ \midrule
    \multirow{6}{*}{\begin{tabular}[c]{@{}c@{}}Switch swap\\ with least\\ insertion\end{tabular}} & $1$ & Load balancing & -7.28 & 56.35 & 0.39 & -0.11 & 2.77 & 565.01 & 125.95 & 64.11 \\
     &  & Rec. Insertion & -8.92 & 7.57 & -0.42 & -1.04 & 4.16 & 627.48 & 111.96 & 60.00 \\ \cmidrule{2-11} 
     & \multirow{2}{*}{$2$} & Load balancing & -7.77 & 5.58 & -0.92 & -0.69 & 3.05 & 850.92 & 128.82 & 64.84 \\
     &  & Rec. Insertion & -10.04 & 9.32 & -0.59 & -0.30 & 3.76 & 823.75 & 132.47 & 57.07 \\ \cmidrule{2-11} 
     & \multirow{2}{*}{$3$} & Load balancing & -6.64 & 5.58 & -1.21 & -0.61 & 2.80 & 1226.21 & 138.13 & 78.00 \\
     &  & Rec. Insertion & -9.51 & 9.20 & -1.01 & -0.73 & 4.96 & 996.88 & 174.79 & 96.07 \\ \midrule
    \multirow{6}{*}{\begin{tabular}[c]{@{}c@{}}Switch swap\\ with least\\ estimated\\ tour\end{tabular}} & \multirow{2}{*}{$1$} & Load balancing & -6.70 & 2.98 & -1.29 & -0.75 & 3.37 & 786.76 & 126.62 & 66.01 \\
     &  & Rec. Insertion & -8.84 & 8.46 & -0.39 & -0.60 & 3.24 & 606.09 & 104.47 & 67.01 \\ \cmidrule{2-11} 
     & \multirow{2}{*}{$2$} & Load balancing & -7.49 & 56.35 & -0.12 & -1.18 & 2.29 & 684.93 & 137.96 & 83.07 \\
     &  & Rec. Insertion & -10.51 & 9.20 & -0.59 & -0.48 & 5.62 & 594.89 & 124.81 & 73.01 \\ \cmidrule{2-11} 
     & \multirow{2}{*}{$3$} & Load balancing & -7.82 & 5.58 & -1.36 & -0.50 & 3.20 & 1171.89 & 204.48 & 104.48 \\
     &  & Rec. Insertion & -10.28 & 9.20 & -1.02 & -1.24 & 3.24 & 1015.93 & 152.57 & 72.04 \\ \midrule
    \multirow{6}{*}{\begin{tabular}[c]{@{}c@{}}Switch swap\\ with least\\ actual tour\end{tabular}} & \multirow{2}{*}{$1$} & Load balancing & -3.99 & 74.20 & 9.90 & 6.53 & 2.55 & 255.50 & 60.99 & 35.37 \\
     &  & Rec. Insertion & -6.46 & 84.85 & 17.39 & 9.28 & 3.04 & 256.34 & 52.93 & 28.78 \\ \cmidrule{2-11} 
     & \multirow{2}{*}{$2$} & Load balancing & -5.14 & 50.19 & 5.29 & 2.16 & 3.03 & 594.80 & 106.44 & 46.70 \\
     &  & Rec. Insertion & -9.27 & 53.03 & 8.39 & 2.70 & 3.51 & 487.24 & 111.47 & 50.54 \\ \cmidrule{2-11} 
     & \multirow{2}{*}{$3$} & Load balancing & -7.51 & 21.39 & 2.22 & 0.97 & 3.29 & 687.75 & 149.73 & 66.34 \\
     &  & Rec. Insertion & -8.29 & 38.33 & 5.00 & 0.67 & 4.54 & 960.34 & 162.15 & 44.35 \\ \bottomrule
    \end{tabular}
\end{sidewaystable}

\begin{sidewaystable}[htb!]
    \centering
    \caption{Comparison of deviation and computation time for three target allocation instances for the variations in the switch and swap neighborhoods}
    \label{tab: three_target_allocation_soln_comparison_switch_swap}
    \begin{tabular}{ccccccccccc}
    \hline
    \multirow{2}{*}{\textbf{Metric}} & \multirow{2}{*}{$n$} & \multirow{2}{*}{\textbf{Construction}} & \multicolumn{4}{c}{\textbf{Deviation from modified MD ($\%$)}} & \multicolumn{4}{c}{\textbf{Computation time (s)}} \\ \cmidrule{4-11} 
     &  &  & Min & Max & Mean & Median & Min & Max & Mean & Median \\ \hline
    \begin{tabular}[c]{@{}c@{}}Modified\\ MD\end{tabular} & - & Load balancing & - & - & - & - & 1.13 & 3281.76 & 676.02 & 73.49 \\ \hline
    \multirow{6}{*}{\begin{tabular}[c]{@{}c@{}}Switch swap\\ with least\\ insertion\end{tabular}} & $1$ & Load balancing & -5.81 & 1.78 & -0.60 & 0.00 & 1.31 & 846.90 & 148.43 & 65.42 \\
     &  & Rec. Insertion & -10.18 & 2.43 & -1.75 & 0.00 & 1.20 & 408.70 & 82.83 & 35.22 \\ \cmidrule{2-11} 
     & \multirow{2}{*}{$2$} & Load balancing & -5.13 & 2.08 & -0.60 & 0.00 & 1.40 & 789.92 & 140.65 & 64.98 \\
     &  & Rec. Insertion & -11.11 & 1.38 & -1.79 & 0.00 & 1.20 & 519.84 & 122.41 & 22.79 \\ \cmidrule{2-11} 
     & \multirow{2}{*}{$3$} & Load balancing & -4.52 & 2.00 & -0.57 & 0.00 & 1.31 & 756.05 & 143.96 & 64.89 \\
     &  & Rec. Insertion & -10.33 & 2.11 & -1.71 & 0.00 & 1.27 & 857.67 & 148.31 & 24.08 \\ \hline
    \multirow{6}{*}{\begin{tabular}[c]{@{}c@{}}Switch swap\\ with least\\ estimated\\ tour\end{tabular}} & \multirow{2}{*}{$1$} & Load balancing & -4.54 & 0.82 & -0.66 & 0.00 & 1.18 & 861.61 & 129.26 & 53.41 \\
     &  & Rec. Insertion & -10.68 & 2.68 & -1.70 & 0.00 & 1.17 & 341.41 & 81.24 & 18.49 \\ \cmidrule{2-11} 
     & \multirow{2}{*}{$2$} & Load balancing & -5.78 & 1.04 & -0.73 & 0.00 & 1.21 & 1014.95 & 188.12 & 50.69 \\
     &  & Rec. Insertion & -10.71 & 1.93 & -1.78 & 0.00 & 1.20 & 847.77 & 126.24 & 17.97 \\ \cmidrule{2-11} 
     & \multirow{2}{*}{$3$} & Load balancing & -5.36 & 1.91 & -0.61 & 0.00 & 1.14 & 1536.28 & 212.81 & 69.64 \\
     &  & Rec. Insertion & -10.77 & 2.33 & -1.96 & -0.09 & 1.20 & 1429.66 & 176.04 & 18.48 \\ \hline
    \multirow{6}{*}{\begin{tabular}[c]{@{}c@{}}Switch swap\\ with least\\ actual tour\end{tabular}} & \multirow{2}{*}{$1$} & Load balancing & -4.14 & 2.45 & -0.26 & 0.00 & 1.07 & 870.21 & 137.39 & 80.54 \\
     &  & Rec. Insertion & -8.19 & 4.91 & -1.25 & 0.00 & 1.18 & 482.32 & 78.53 & 35.35 \\ \cmidrule{2-11} 
     & \multirow{2}{*}{$2$} & Load balancing & -6.32 & 2.62 & -0.37 & 0.00 & 1.07 & 651.75 & 137.15 & 73.31 \\
     &  & Rec. Insertion & -9.95 & 2.12 & -1.62 & 0.00 & 1.17 & 521.69 & 106.90 & 43.73 \\ \cmidrule{2-11} 
     & \multirow{2}{*}{$3$} & Load balancing & -6.10 & 3.53 & -0.57 & 0.00 & 1.12 & 1481.60 & 233.39 & 52.02 \\
     &  & Rec. Insertion & -10.94 & 5.07 & -1.70 & 0.00 & 1.17 & 1172.43 & 157.72 & 18.49 \\ \hline
    \end{tabular}
\end{sidewaystable}

\begin{sidewaystable}[htb!]
    \centering
    \caption{Comparison of deviation and computation time for five target allocation instances for the variations in the switch and swap neighborhoods}
    \label{tab: five_target_allocation_soln_comparison_switch_swap}
    \begin{tabular}{ccccccccccc}
    \hline
    \multirow{2}{*}{\textbf{Metric}} & \multirow{2}{*}{$n$} & \multirow{2}{*}{\textbf{Construction}} & \multicolumn{4}{c}{\textbf{Deviation from modified MD ($\%$)}} & \multicolumn{4}{c}{\textbf{Computation time (s)}} \\ \cmidrule{4-11} 
     &  &  & Min & Max & Mean & Median & Min & Max & Mean & Median \\ \hline
    \begin{tabular}[c]{@{}c@{}}Modified\\ MD\end{tabular} & - & Load balancing & - & - & - & - & 1.29 & 3600.00 & 681.84 & 104.14 \\ \hline
    \multirow{6}{*}{\begin{tabular}[c]{@{}c@{}}Switch swap\\ with least\\ insertion\end{tabular}} & $1$ & Load balancing & -4.90 & 1.71 & -0.59 & -0.01 & 1.53 & 606.92 & 138.64 & 60.47 \\
     &  & Rec. Insertion & -8.65 & 0.63 & -1.59 & 0.00 & 1.35 & 453.47 & 82.23 & 19.38 \\ \cmidrule{2-11} 
     & \multirow{2}{*}{$2$} & Load balancing & -4.03 & 1.10 & -0.64 & 0.00 & 1.50 & 522.75 & 138.41 & 51.78 \\
     &  & Rec. Insertion & -8.45 & 0.56 & -1.55 & -0.01 & 1.31 & 439.40 & 105.12 & 16.18 \\ \cmidrule{2-11} 
     & \multirow{2}{*}{$3$} & Load balancing & -4.41 & 0.63 & -0.47 & 0.00 & 1.51 & 587.57 & 137.00 & 61.56 \\
     &  & Rec. Insertion & -8.62 & 2.47 & -1.41 & 0.00 & 1.23 & 911.96 & 125.11 & 19.80 \\ \hline
    \multirow{6}{*}{\begin{tabular}[c]{@{}c@{}}Switch swap\\ with least\\ estimated\\ tour\end{tabular}} & \multirow{2}{*}{$1$} & Load balancing & -4.06 & 3.42 & -0.45 & 0.00 & 1.88 & 579.03 & 113.24 & 46.34 \\
     &  & Rec. Insertion & -8.11 & 1.20 & -1.48 & 0.00 & 1.22 & 452.33 & 82.44 & 17.67 \\ \cmidrule{2-11} 
     & \multirow{2}{*}{$2$} & Load balancing & -3.33 & 4.04 & -0.36 & 0.00 & 1.31 & 734.49 & 155.91 & 37.89 \\
     &  & Rec. Insertion & -8.05 & 2.66 & -1.44 & 0.00 & 1.35 & 850.86 & 114.09 & 18.36 \\ \cmidrule{2-11} 
     & \multirow{2}{*}{$3$} & Load balancing & -5.52 & 4.64 & -0.38 & 0.00 & 1.64 & 918.58 & 175.57 & 35.54 \\
     &  & Rec. Insertion & -8.41 & 0.04 & -1.52 & -0.01 & 1.24 & 610.09 & 117.44 & 18.77 \\ \hline
    \multirow{6}{*}{\begin{tabular}[c]{@{}c@{}}Switch swap\\ with least\\ actual tour\end{tabular}} & \multirow{2}{*}{$1$} & Load balancing & -2.45 & 18.03 & 0.69 & 0.00 & 1.38 & 651.61 & 113.74 & 34.69 \\
     &  & Rec. Insertion & -8.26 & 2.61 & -1.14 & 0.00 & 1.24 & 314.71 & 68.64 & 16.32 \\ \cmidrule{2-11} 
     & \multirow{2}{*}{$2$} & Load balancing & -3.86 & 6.58 & -0.17 & 0.00 & 1.62 & 1087.70 & 148.03 & 45.14 \\
     &  & Rec. Insertion & -8.35 & 2.66 & -1.46 & 0.00 & 1.21 & 601.64 & 107.08 & 18.69 \\ \cmidrule{2-11} 
     & \multirow{2}{*}{$3$} & Load balancing & -3.86 & 4.67 & -0.33 & 0.00 & 1.76 & 1003.26 & 189.54 & 49.03 \\
     &  & Rec. Insertion & -8.18 & 1.26 & -1.56 & -0.01 & 1.20 & 794.78 & 135.69 & 18.99 \\ \hline
    \end{tabular}
\end{sidewaystable}

\subsection{Multi-target swap neighborhood analysis} \label{appsubsect: comp_results_mult_target}

The detailed results for the multi-target swap neighborhood analysis described in Section~\ref{subsect: mult_targ_swap} are provided in this section. In particular, the description of the minimum, maximum, mean, and median deviation of the objective value with respect to the modified MD algorithm and the minimum, maximum, mean, and median computation time are provided in Tables~\ref{tab: no_target_allocation_soln_comparison_mult_targ}, \ref{tab: three_target_allocation_soln_comparison_mult_targ}, and \ref{tab: five_target_allocation_soln_comparison_mult_targ}, for the zero-target, three-target, and five-target instances, respectively, selected based on the sensitivity study in Section~\ref{subsect: sensitive_instances}.

\begin{sidewaystable}[htb!]
    \centering
    \caption{Comparison of deviation and computation time for zero-target allocation instances for the variations in the multi-target swap neighborhood}
    \label{tab: no_target_allocation_soln_comparison_mult_targ}
    \begin{tabular}{ccccccccccc}
    \toprule
    \multirow{2}{*}{\textbf{Metric}} & \multirow{2}{*}{$m$} & \multirow{2}{*}{$n$} & \multicolumn{4}{c}{\textbf{Deviation from modified MD (\%)}} & \multicolumn{4}{c}{\textbf{Computation time (s)}} \\ \cmidrule{4-11} 
     &  &  & Min & Max & Mean & Median & Min & Max & Mean & Median \\ \midrule
    Modified MD & - & - & - & - & - & - & 30.59 & 3415.71 & 1140.64 & 893.61 \\ \midrule
    \begin{tabular}[c]{@{}c@{}}Switch swap with \\ least insertion\end{tabular} & - & - & -7.77 & 4.10 & -0.93 & -0.32 & 8.08 & 376.60 & 171.84 & 162.19 \\ \midrule
    \multirow{4}{*}{\begin{tabular}[c]{@{}c@{}}Switch swap with fixed \\ struct. mult. targ. swap\\ using insertion cost\end{tabular}} & \multirow{2}{*}{$2$} & $10$ & -10.07 & 7.47 & -0.29 & -0.44 & 5.91 & 1134.32 & 403.07 & 283.38 \\
     &  & $20$ & -9.83 & 9.21 & -0.75 & -0.75 & 5.71 & 1183.31 & 460.81 & 354.44 \\ \cmidrule{2-11} 
     & \multirow{2}{*}{$2, 3$} & $10$ & -9.85 & 6.87 & -1.05 & -1.62 & 7.44 & 2827.50 & 855.40 & 480.28 \\
     &  & $20$ & -9.51 & 9.63 & -1.20 & -1.18 & 6.40 & 3172.68 & 886.35 & 502.12 \\ \midrule
    \multirow{4}{*}{\begin{tabular}[c]{@{}c@{}}Switch swap with fixed\\ struct. mult. targ. swap\\ using savings and\\ insertion cost\end{tabular}} & \multirow{2}{*}{$2$} & $10$ & -9.58 & 7.82 & -0.55 & -0.02 & 6.62 & 1519.80 & 584.24 & 395.95 \\
     &  & $20$ & -10.51 & 9.08 & -0.51 & -0.30 & 7.34 & 1419.80 & 516.78 & 321.38 \\ \cmidrule{2-11} 
     & \multirow{2}{*}{$2, 3$} & $10$ & -9.89 & 5.78 & -1.65 & -1.67 & 7.19 & 2521.33 & 889.43 & 433.93 \\
     &  & $20$ & -11.10 & 9.20 & -1.87 & -1.48 & 11.06 & 3600.00 & 1238.71 & 682.95 \\ \midrule
    \multirow{4}{*}{\begin{tabular}[c]{@{}c@{}}Switch swap with var.\\ struct. mult. targ. swap\\ using saving and insertion\\ cost and recur. insertion\end{tabular}} & \multirow{2}{*}{$3$} & $10$ & -10.66 & 7.82 & -1.79 & -1.89 & 6.38 & 2618.62 & 993.83 & 737.75 \\
     &  & $20$ & -10.19 & 4.27 & -2.22 & -1.92 & 7.02 & 3600.00 & 1414.77 & 665.21 \\ \cmidrule{2-11} 
     & \multirow{2}{*}{$4$} & $10$ & -10.80 & 6.63 & -0.75 & -0.99 & 7.12 & 3600.00 & 1896.57 & 1266.91 \\
     &  & $20$ & -11.25 & 9.20 & -1.51 & -0.86 & 6.95 & 3600.00 & 1852.79 & 1680.34 \\ \midrule
    \multirow{4}{*}{\begin{tabular}[c]{@{}c@{}}Switch swap with var.\\ struct. mult. targ. swap\\ using saving and insertion\\ cost and fixed insertion\end{tabular}} & \multirow{2}{*}{$3$} & $10$ & -10.30 & 9.20 & -1.08 & -1.47 & 7.67 & 2714.96 & 926.48 & 735.36 \\
     &  & $20$ & -10.44 & 8.73 & -1.59 & -1.93 & 4.51 & 3637.48 & 1253.17 & 536.09 \\ \cmidrule{2-11} 
     & \multirow{2}{*}{$4$} & $10$ & -11.30 & 6.01 & -1.15 & -2.48 & 6.81 & 3600.00 & 1430.14 & 944.76 \\
     &  & $20$ & -10.88 & 9.20 & -1.66 & -1.94 & 6.40 & 3646.12 & 1649.91 & 1136.80 \\ \bottomrule
    \end{tabular}
\end{sidewaystable}

\begin{sidewaystable}[htb!]
    \centering
    \caption{Comparison of deviation and computation time for three-target allocation instances for the variations in the multi-target swap neighborhood}
    \label{tab: three_target_allocation_soln_comparison_mult_targ}
    \begin{tabular}{ccccccccccc}
    \toprule
    \multirow{2}{*}{\textbf{Metric}} & \multirow{2}{*}{$m$} & \multirow{2}{*}{$n$} & \multicolumn{4}{c}{\textbf{Deviation from modified MD (\%)}} & \multicolumn{4}{c}{\textbf{Computation time (s)}} \\ \cmidrule{4-11} 
     &  &  & Min & Max & Mean & Median & Min & Max & Mean & Median \\ \midrule
    Modified MD & - & - & - & - & - & - & 5.08 & 3050.48 & 1151.84 & 1237.26 \\ \midrule
    \begin{tabular}[c]{@{}c@{}}Switch swap with \\ least insertion\end{tabular} & - & - & -3.94 & 0.84 & -0.74 & -0.49 & 6.11 & 622.04 & 225.33 & 190.77 \\ \midrule
    \multirow{4}{*}{\begin{tabular}[c]{@{}c@{}}Switch swap with fixed \\ struct. mult. targ. swap\\ using insertion cost\end{tabular}} & \multirow{2}{*}{$2$} & $10$ & -5.87 & 2.07 & -2.26 & -2.67 & 5.16 & 1015.93 & 405.35 & 374.63 \\
     &  & $20$ & -6.22 & 1.30 & -2.71 & -2.16 & 4.06 & 1546.86 & 523.55 & 336.73 \\ \cmidrule{2-11} 
     & \multirow{2}{*}{$2, 3$} & $10$ & -5.50 & 1.90 & -2.22 & -2.33 & 5.24 & 1932.48 & 561.15 & 326.62 \\
     &  & $20$ & -5.59 & 0.72 & -2.30 & -2.29 & 4.58 & 2620.95 & 705.99 & 374.12 \\ \midrule
    \multirow{4}{*}{\begin{tabular}[c]{@{}c@{}}Switch swap with fixed\\ struct. mult. targ. swap\\ using savings and\\ insertion cost\end{tabular}} & \multirow{2}{*}{$2$} & $10$ & -6.22 & 2.16 & -2.08 & -2.51 & 4.54 & 1250.50 & 463.65 & 290.29 \\
     &  & $20$ & -6.22 & 0.52 & -2.34 & -2.05 & 5.30 & 1236.33 & 458.71 & 330.90 \\ \cmidrule{2-11} 
     & \multirow{2}{*}{$2, 3$} & $10$ & -5.58 & 1.91 & -2.05 & -2.10 & 3.38 & 1655.43 & 593.14 & 356.83 \\
     &  & $20$ & -6.25 & 1.78 & -2.28 & -2.04 & 4.26 & 2721.16 & 842.91 & 428.35 \\ \midrule
    \multirow{4}{*}{\begin{tabular}[c]{@{}c@{}}Switch swap with var.\\ struct. mult. targ. swap\\ using saving and insertion\\ cost and recur. insertion\end{tabular}} & \multirow{2}{*}{$3$} & $10$ & -6.46 & 1.68 & -2.64 & -2.78 & 5.77 & 3600.00 & 1336.91 & 1019.23 \\
     &  & $20$ & -6.05 & 0.93 & -2.27 & -2.15 & 3.56 & 1735.64 & 848.99 & 755.53 \\ \cmidrule{2-11} 
     & \multirow{2}{*}{$4$} & $10$ & -6.27 & 0.63 & -2.71 & -2.47 & 6.06 & 3600.00 & 2072.46 & 2441.88 \\
     &  & $20$ & -6.02 & 2.12 & -2.18 & -2.38 & 2.25 & 3600.00 & 1628.39 & 1529.93 \\ \midrule
    \multirow{4}{*}{\begin{tabular}[c]{@{}c@{}}Switch swap with var.\\ struct. mult. targ. swap\\ using saving and insertion\\ cost and fixed insertion\end{tabular}} & \multirow{2}{*}{$3$} & $10$ & -6.15 & 1.48 & -2.32 & -2.17 & 5.03 & 2159.17 & 772.92 & 578.04 \\
     &  & $20$ & -6.11 & 1.52 & -2.56 & -2.53 & 5.27 & 2305.64 & 971.47 & 563.82 \\ \cmidrule{2-11} 
     & \multirow{2}{*}{$4$} & $10$ & -6.36 & 1.13 & -2.45 & -2.47 & 4.22 & 3600.00 & 1619.27 & 1185.42 \\
     &  & $20$ & -5.64 & 1.28 & -2.49 & -2.57 & 3.73 & 3600.00 & 1592.65 & 1294.97 \\ \bottomrule
    \end{tabular}
\end{sidewaystable}

\begin{sidewaystable}[htb!]
    \centering
    \caption{Comparison of deviation and computation time for five-target allocation instances for the variations in the multi-target swap neighborhood}
    \label{tab: five_target_allocation_soln_comparison_mult_targ}
    \begin{tabular}{ccccccccccc}
    \toprule
    \multirow{2}{*}{\textbf{Metric}} & \multirow{2}{*}{$m$} & \multirow{2}{*}{$n$} & \multicolumn{4}{c}{\textbf{Deviation from modified MD (\%)}} & \multicolumn{4}{c}{\textbf{Computation time (s)}} \\ \cmidrule{4-11} 
     &  &  & Min & Max & Mean & Median & Min & Max & Mean & Median \\ \midrule
    Modified MD & - & - & - & - & - & - & 1.29 & 3489.23 & 1309.79 & 895.87 \\ \midrule
    \begin{tabular}[c]{@{}c@{}}Switch swap with \\ least insertion\end{tabular} & - & - & -4.03 & 0.92 & -1.21 & -1.22 & 1.50 & 522.75 & 242.00 & 251.55 \\ \midrule
    \multirow{4}{*}{\begin{tabular}[c]{@{}c@{}}Switch swap with fixed \\ struct. mult. targ. swap\\ using insertion cost\end{tabular}} & \multirow{2}{*}{$2$} & $10$ & -7.42 & 1.03 & -2.59 & -2.53 & 1.25 & 1570.85 & 391.88 & 252.50 \\
     &  & $20$ & -6.84 & 1.32 & -2.72 & -2.49 & 1.22 & 2252.52 & 475.26 & 266.62 \\ \cmidrule{2-11} 
     & \multirow{2}{*}{$2, 3$} & $10$ & -7.08 & 1.12 & -2.82 & -2.26 & 1.25 & 1968.77 & 478.96 & 340.14 \\
     &  & $20$ & -7.10 & 0.70 & -2.34 & -1.63 & 1.32 & 1737.77 & 415.27 & 267.66 \\ \midrule
    \multirow{4}{*}{\begin{tabular}[c]{@{}c@{}}Switch swap with fixed\\ struct. mult. targ. swap\\ using savings and\\ insertion cost\end{tabular}} & \multirow{2}{*}{$2$} & $10$ & -7.45 & 0.61 & -2.71 & -2.20 & 1.24 & 1186.64 & 387.84 & 349.06 \\
     &  & $20$ & -6.72 & 1.54 & -2.62 & -2.27 & 1.31 & 1356.44 & 351.47 & 183.38 \\ \cmidrule{2-11} 
     & \multirow{2}{*}{$2, 3$} & $10$ & -7.21 & 0.91 & -2.64 & -2.16 & 1.31 & 3281.01 & 651.16 & 327.88 \\
     &  & $20$ & -7.12 & 0.88 & -2.75 & -2.41 & 1.29 & 2598.03 & 583.66 & 286.08 \\ \midrule
    \multirow{4}{*}{\begin{tabular}[c]{@{}c@{}}Switch swap with var.\\ struct. mult. targ. swap\\ using saving and insertion\\ cost and recur. insertion\end{tabular}} & \multirow{2}{*}{$3$} & $10$ & -7.03 & 0.00 & -3.15 & -2.17 & 1.30 & 3600.00 & 988.02 & 524.34 \\
     &  & $20$ & -7.98 & 0.56 & -2.90 & -2.43 & 1.33 & 3600.00 & 983.70 & 624.23 \\ \cmidrule{2-11} 
     & \multirow{2}{*}{$4$} & $10$ & -6.91 & 0.54 & -3.05 & -2.63 & 1.39 & 3600.00 & 1501.73 & 1108.68 \\
     &  & $20$ & -6.87 & 0.00 & -3.23 & -3.02 & 1.32 & 3600.00 & 1316.45 & 945.60 \\ \midrule
    \multirow{4}{*}{\begin{tabular}[c]{@{}c@{}}Switch swap with var.\\ struct. mult. targ. swap\\ using saving and insertion\\ cost and fixed insertion\end{tabular}} & \multirow{2}{*}{$3$} & $10$ & -6.77 & 0.67 & -3.05 & -2.65 & 1.30 & 2881.29 & 717.13 & 449.93 \\
     &  & $20$ & -7.14 & 0.87 & -2.82 & -2.40 & 1.37 & 3499.93 & 838.21 & 406.42 \\ \cmidrule{2-11} 
     & \multirow{2}{*}{$4$} & $10$ & -7.52 & 0.89 & -3.13 & -2.56 & 1.23 & 3600.00 & 1267.56 & 761.42 \\
     &  & $20$ & -6.95 & 0.70 & -3.10 & -2.70 & 1.34 & 3600.00 & 956.48 & 663.56 \\ \bottomrule
    \end{tabular}
\end{sidewaystable}

\end{appendices}

\bibliography{cite}

\end{document}